\documentclass[preprint,12pt]{article}

\usepackage[english]{babel}	
\usepackage[utf8x]{inputenc}		
\usepackage{enumerate}
\usepackage{textcomp}                
\usepackage{typearea}
\usepackage{pdfpages}
\usepackage[toc]{appendix}
\usepackage[]{todonotes}
\usepackage{subcaption}

\usepackage{amsfonts, amsmath, amssymb, amsthm, dsfont}
\usepackage{amsthm}			
\usepackage{thumbpdf}	
\usepackage{booktabs}

\usepackage{pgf, tikz}
\usepackage{tikz-cd}		
\usepackage{bbm}

\theoremstyle{definition}
\newtheorem{definition}{Definition}[section]

\newtheorem*{assumption*}{Assumption}

\newtheorem*{condition*}{Condition}

\theoremstyle{plain}
\newtheorem{theorem}[definition]{Theorem}
\newtheorem{proposition}[definition]{Proposition}
\newtheorem{lemma}[definition]{Lemma}
\newtheorem{cor}[definition]{Corollary}
\theoremstyle{remark}
\newtheorem{remark}{Remark}

\usepackage{graphicx}
\usepackage{hyperref}

\newcommand{\N}{\mathbb{N}}

\newcommand{\Q}{\mathbb{Q}}
\newcommand{\R}{\mathbb{R}}

\newcommand{\E}{\mathbb{E}}

\newcommand{\F}{\mathcal{F}}

\newcommand{\pconv}{\xrightarrow{P}}

\newcommand{\deq}{\overset{d}{=}}

\newcommand{\wconv}{\Rightarrow}

\newcommand{\sgn}{\mathrm{sign}}

\newcommand{\vs}{\mathbf{s}}

\newcommand{\vt}{\mathbf{t}}
\newcommand{\vr}{\mathbf{r}}

\title{Regularity of multifractional moving average processes with random Hurst exponent\footnote{Accepted manuscript, to appear in \textit{Stochastic Processes and their Applications}, DOI:10.1016/j.spa.2021.05.008}}
\author{Dennis Loboda\thanks{loboda@stochastik.rwth-aachen.de} \qquad Fabian Mies\thanks{mies@stochastik.rwth-aachen.de (corresponding author)} \qquad Ansgar Steland\thanks{steland@stochastik.rwth-aachen.de} \\[2ex] RWTH Aachen University\\ Institute of Statistics}

\begin{document}

\maketitle

\begin{abstract}
	A recently proposed alternative to multifractional Brownian motion (mBm) with random Hurst exponent is studied, which we refer to as Itô-mBm.
	It is shown that Itô-mBm is locally self-similar.
	In contrast to mBm, its pathwise regularity is almost unaffected by the roughness of the functional Hurst parameter.
	The pathwise properties are established via a new polynomial moment condition similar to the Kolmogorov-Centsov theorem, allowing for random local Hölder exponents.	
	Our results are applicable to a broad class of moving average processes where pathwise regularity and long memory properties may be decoupled, e.g.\ to a multifractional generalization of the Matérn process.
	
\end{abstract}

\section{Introduction}

Regularity or roughness of the sample paths of a stochastic process $ Y_t, t \ge 0, $ is a crucial property, e.g., when studying models defined by stochastic integrals with respect to (w.r.t.) $ Y_t $. It can be measured in terms of the (local) $p$-variation, a non-increasing function of $ 0 < p < \infty $, and paths with infinite $1$-variation are called rough. Properties and the related calculus of rough functions is reviewed and discussed in \cite{Norvaisa2006}. Alternatively, one may establish Hölder-continuity and make use of the fact that $\alpha$-Hölder continuity implies bounded $1/\alpha$-variation. A convenient criterion is the Kolmogorov-Centsov Theorem (see e.g.\ \cite[Thm.\ 2.8]{Karatzas1998}), which links Hölder-regularity of the paths to absolute moments of the increments. This well-known result can be used to establish Hölder coefficients which (i) hold uniformly on the interval of interest and (ii) are deterministic. 
This is satisfactory for most conventional processes of interest, e.g.\ fractional Brownian motion or solutions of stochastic differential equations. 
For general stochastic processes, the regularity of sample paths neither needs to be deterministic nor globally constant.
A very simple example is the stopped Brownian motion $(W_{t\wedge \tau})_{t\geq 0}$ for a random stopping time $\tau>0$. 
On $[0,\tau)$, the process is (almost) $1/2$ Hölder continuous, whereas it is arbitrarily smooth on $(\tau,\infty)$ as a constant.
Hence, the smoothness at any fixed time $t$ is random.

The  process primarily considered in this paper is closely related to the class of multifractional Brownian motion (mBm), which generalizes fractional Brownian motion (fBm) by allowing the Hurst exponent to vary in time. Multifractional Brownian motions have been used as a model for stochastic volatility in finance \cite{Corlay2014}, for network traffic \cite{Bianchi2004}, and for temperature time series \cite{LevyVehel2013}, among others.
As there exist various definitions of mBm in the literature leading to slightly different processes \cite{Peltier1995,Benassi1997,Stoev2006,Lebovits2014}, we shall confine our discussion to a special case of \cite[Def.\ 1.1]{Stoev2006}.
Let $W_t$ be a Brownian motion defined on a filtered probability space $(\Omega, (\F_t)_{t\in\R}, \F, P)$, and define \begin{align}
	B(t,H) = \int_{-\infty}^{t} (t-s)_+^{H-\frac{1}{2}}-(-s)_+^{H-\frac{1}{2}}\, dW_s,\quad t\in\R, H\in(0,1), \label{eqn:def-fbm}
\end{align}
which is a random field indexed by $t$ and $H$. For fixed $H$, $B(t,H)$ is a fBm with parameter $H$.
For a deterministic process $H_t$, the mBm may be defined as $B^H_t = B(t, H_t)$, i.e.
\begin{align}
	B_t^H = B(t, H_t) = \int_{-\infty}^{t} (t-s)_+^{H_t-\frac{1}{2}}-(-s)_+^{H_t-\frac{1}{2}}\, dW_s,\quad t\in\R. \label{eqn:mbm-field}
\end{align}
Intuitively, $B_t^H$ behaves like a fractional Brownian motion with Hurst parameter $H_t$ locally around $t$.

To describe the local regularity of a stochastic process such as $B^H_t$, we define for a generic process $Y_t$ the pointwise Hölder coefficient as
\begin{align}
	\alpha_t(Y) = \sup\left\{ \alpha \, :\,  \limsup_{\epsilon\to 0} \sup_{|h|<\epsilon} \frac{|Y_{t+h}-Y_t|}{|h|^\alpha} = 0 \right\}. \label{eqn:def-pointwise-exponent}
\end{align}
It has been shown by \cite[Prop.\ 13]{Herbin2006} that for any $t$, almost surely, \begin{align}
	\alpha_t(B^H) = H_t \wedge \alpha_t(H).\label{eqn:mbm-pointwise-smoothness}
\end{align}
That is, $\alpha_t(B^H)$ is a modification of $H_t\wedge \alpha_t(H)$, but the latter two processes are in general not indistinguishable \cite{Ayache2013}.

The identity \eqref{eqn:mbm-pointwise-smoothness} for the pointwise Hölder exponent reveals that the regularity of the functional Hurst exponent $H_t$ itself significantly affects the smoothness of $B^H_t$.
In the extreme case that $H_t$ is a discontinuous step-function, we readily find from \eqref{eqn:mbm-field} that $B^H_t$ itself is discontinuous.
In particular, $t\mapsto B(t,H)$ is continuous for all $H\in(0,1)$, but $B(t,H)\neq B(t, H')$ for $H\neq H'$ because the integrands differ. 
This discontinuity invalidates the intuition that $B^H_t$ should behave locally like a fBm with Hurst parameter $H_t$.
The relevance of the regularity of $H_t$ even extends to the statistical estimation of $H_t$ for a fixed $t$.
For example, the rates of convergence of the nonparametric estimators of \cite{Bardet2013} are vacuous unless $\alpha_t(H)>H_t$ (see Proposition 3 therein).
While nonparametric estimators always require some smoothness of the target quantity, i.e.\ $\alpha_t(H)$, the relevance of the difference $\alpha_t(H)-H_t$ is a special feature of inference for mBm.
When estimating the global Hölder regularity $\inf_t H_t$, \cite{Lebovits2017} also require that $H_t$ is $\eta$-Hölder continuous for some $\eta>\sup_t H_t$.

Another shortcoming of defining a multifractional process in terms of \eqref{eqn:mbm-field} is that this construction can not easily be extended to a stochastic ($\F_t$-adapted) process $H_t$. 
In fact, if $H_t$ in \eqref{eqn:mbm-field} is random, the integrand is in general no longer $\F_s$-adapted and the stochastic integral can not be interpreted in the Itô sense.
A possible remedy is to either assume $H_t$ to be independent of the driving Gaussian noise $(W_t)_{t\in\R}$, or to employ an alternative representation of the random field $B(t,H)$ \cite{Ayache2005}.

As an alternative to overcome the latter measurability issues, \cite{Ayache2018} suggest to instead study the process\footnote{In \cite{Ayache2018}, only the term $\int_0^t (t-s)_+^{H_s-1/2}\, dW_s$ is considered in detail.}
\begin{align}
	K^H_t = \int_{-\infty}^{t} \sigma_s\left[(t-s)_+^{H_s-\frac{1}{2}} - (-s)_+^{H_s-\frac{1}{2}}\right]\, dW_s. \label{eqn:mbm-ito}
\end{align}
By making the integrand depend on $H_s$ instead of $H_t$, the integral in \eqref{eqn:mbm-ito} can be defined in the classical Itô sense.
A similar process has also been proposed by \cite{Surgailis2008}, though for a deterministic Hurst function $H_s$.
The process $K_t^H$ has been called a multifractional process with random exponent (MPRE) by \cite{Ayache2018}.
To contrast this process with mBm, we will refer to $K_t^H$ as an Itô-mBm.
To study the regularity of $K^H_t$, \cite{Ayache2018} derive a series representation based on a wavelet basis, which can be investigated via analytical methods. 
Their findings suggest that the smoothness of \eqref{eqn:mbm-ito} is less sensitive to the regularity of $H_t$.
In particular, under the condition that $\inf_tH_t>1/2$ and that $H_t$ is $\eta$-Hölder continuous for $\eta>1/2$, they find that the Hölder regularity of $K^H_t$ on any interval $[t_1, t_2]$ is at least $\inf_{s\in[t_1, t_2]} H_s$. 
Although the applicability of this result is restricted, it allows for the case $\eta<H_t$, much unlike \eqref{eqn:mbm-pointwise-smoothness}.

The regularity results derived by \cite{Ayache2018} are based on a wavelet representation of $K_t^H$. 
In contrast, the novel Kolmogorov-Centsov type result presented in Section \ref{sec:KC} allows us to simplify the analysis of the Itô-mBm by probabilistic methods, and we are able to obtain refined results on its pathwise regularity.
In particular, we show in Section \ref{sec:mbm} that its pointwise Hölder exponent satisfies $\alpha_t(K^H) \geq H_t$ if $H_t$ is continuous, irrespective of the regularity of $H_t$.
Different from the analysis of \cite{Ayache2018}, we may include the cases $H_t\leq \tfrac{1}{2}$ and $\alpha_t(H)\leq \tfrac{1}{2}$.
In contrast to the classical mBm $B^H_t$, the paths of $H_t$ may be rougher than those of the corresponding Itô-mBm $K^H_t$.
For the case that $H_t$ is discontinuous, we also obtain similar results which are only slightly weaker.
Notably, in the latter case, $K^H_t$ is still continuous, in contrast to the process $B_t^H$.

To show that our estimates of the pointwise Hölder exponent are sharp, i.e.\ that $\alpha_t(K^H)=H_t$, we establish a rescaling limit showing that, for each fixed $t$, as $h\to 0$, 
\begin{align}
	h^{-H_t} (K^H_{t+hr}-K^H_t) \; \wconv\;  \sigma_t \int_{-\infty}^{r}(r-s)_+^{H_t-\frac{1}{2}} - (-s)_+^{H_t-\frac{1}{2}}\, d\tilde{W}_s, \label{eqn:rescaling-intro}
\end{align}
for a Brownian motion $\tilde{W}_s$ independent of $H_t$. 
That is, the process $K^H_t$ behaves locally near $t$ like a fractional Brownian motion with Hurst parameter $H_t$.
The precise limit theorem is presented in Section \ref{sec:mbm}.
All our results hold for a broader class of moving average processes which contains, for example, a multifractional generalization of the Mátern process, see Section \ref{sec:mbm}.

The processes $K^H_t$ and $B^H_t$ are in general different, except for the trivial case where $H_t$ is constant. 
Nevertheless, both processes may serve as a non-stationary generalization of fractional Brownian motion. 
In particular, the scaling relation \eqref{eqn:rescaling-intro} also holds for $B_t^{H}$ if $H_t$ is sufficiently smooth, since $H\mapsto B(t,H)$ is $C^\infty$ \cite{Ayache2013}. 
Thus, the presented smoothness results for $K_t^H$ could be an argument to use the latter process in practice.

This paper is structured as follows. 
In Section \ref{sec:KC}, we present the continuity criterion result for random fields with random, local Hölder exponents. 
These results are applied in Section \ref{sec:mbm} to study the smoothness of the Itô-mBm $K^H_t$.
All proofs are gathered in Section \ref{sec:proofs}.

\subsection*{Notation}
For two real numbers $a,b$, we denote $a\wedge b =\min(a,b)$, $a\vee b=\max(a,b)$, and $(a)_+ = \max(a,0)$.
To make clear where we study random fields, we denote scalar indices $s,t \in \R$ by normal letters and vector indices $\vs,\vt \in \R^d$ by bold letters.
We denote by $C$ a generic constant, the value of which might change from line to line.
If the factor $C$ depends on the quantities $a,b,c$, we denote this as $C(a,b,c)$.
Weak convergence in metric spaces is denoted as $\wconv$.
For an interval $I\subset \R$, the space of continuous functions equipped with the supremum norm is denoted by $C(I)$. 
Convergence in probability is denoted as $\pconv$, and for a sequence of random variables $\varepsilon_n$, we write $\varepsilon_n=o_P(1)$ if $\varepsilon_n\pconv 0$ as $n\to\infty$.

\section{Pathwise Hölder continuity of random fields}\label{sec:KC}

On a probability space $(\Omega,\F, P)$, we consider a real-valued random field $(Y_{\vt})_{\vt}$ indexed by $\vt \in[0,T]^d$, i.e.\ each $Y_{\vt}$ is a random variable.
The Kolmogorov-Centsov theorem requires that there are two positive real numbers $\alpha,\beta>0$, and a constant $C$, such that for all $\vs,\vt\in[0,T]^d$, it holds 
\begin{align}
	\E|Y_{\vs}-Y_{\vt}|^\alpha \leq C \|\vs-\vt\|^{d+\beta}, \label{eqn:KC-condition}
\end{align}
where $\|\cdot\|$ is an arbitrary norm on $\R^d$. 
If \eqref{eqn:KC-condition} holds, there exists a modification $\tilde{Y}_{\vt}$ of $Y_{\vt}$, i.e.\ $P(\tilde{Y}_{\vt}=Y_{\vt})=1$, such that the paths of $\tilde{Y}_{\vt}$ are $\eta$-Hölder continuous for any $\eta\in(0, \beta/\alpha)$. 
While this is a result on the global Hölder continuity of $\tilde{Y}_{\vt}$, it can be localized by requiring \eqref{eqn:KC-condition} to hold in a neighborhood of some point ${\vt}_0$.
However, this moment condition does not generalize to random Hölder exponents as it is only sensible for deterministic values $\alpha$ and $\beta$.
As a remedy, we suggest to let $\beta$ be random, and to treat it jointly with the expectation. 
This leads to the moment criterion \eqref{eqn:moment-condition} below.
Note that Theorem \ref{thm:KC} allows for both, randomness and local variation of the Hölder exponent  $a_\vt$, which is thus represented as a random field.
As a regularity assumption, we require $a_\vt$ to be lower semicontinuous, that is
\begin{align}
	\lim_{\epsilon\to 0}\inf_{\substack{\vs\in[0,T]^d \\ \|\vs-\vt\|\leq \epsilon}} a_\vs \quad\geq\quad a_\vt,\qquad \forall\, \vt \in [0,T]^d.\label{eqn:def-LSC}
\end{align}
Furthermore, we introduce the notation $[\vs,\vt]=\{ \mathbf{r} \in [0,T]^d: (s_i\wedge t_i) \leq r_i \leq (s_i\vee t_i) \}$ for the hyper-rectangles induced by $\vt,\vs\in [0,T]^d$.

\begin{theorem}[Local continuity criterion]\label{thm:KC}
	Let $Y_{\vt}\in\R,a_{\vt}\in(0,1), \vt\in [0,T]^d$ be random fields.
	Let $a_\vt$ be such that $\underline{a}_B = \inf_{\vt\in B} a_\vt$ is measurable for any open or closed set $B\subset[0,T]^d$, and assume that $\inf_{\vt \in [0,T]^d} a_{\vt} >0$, and that $a_\vt$ is lower semicontinuous. 
	Suppose that for some $p>0$, there exists an $\epsilon>0$ and a constant $C(p,\epsilon)$ such that for all $\vs,\vt \in[0,T]^d$ with $\|\vs- \vt\|\leq \epsilon$, it holds 
	\begin{align}
		\E  \left| \frac{Y_{\vs}-Y_{\vt}}{\|\vs-\vt \|^{\underline{a}_{[\vs,\vt]}}} \right|^p \leq C(p,\epsilon) \|\vs-\vt\|^d. \label{eqn:moment-condition}
	\end{align}
	Then there exists a modification $\tilde{Y}_{\vt}$ of $Y_{\vt}$ such that, 
	\begin{align}
		P\left( \forall \gamma>0\; \forall B\in\mathcal{B}:  \sup_{\substack{\vs,\vt\in B \\ \vs \neq \vt}} \frac{|\tilde{Y}_{\vs} - \tilde{Y}_{\vt}|}{ \|\vs-\vt\|^{\underline{a}_B-\gamma}} <\infty \right) = 1, \label{eqn:KC-result}
	\end{align}
	where $\mathcal{B}$ is the collection of closed subsets of $[0,T]^d$.
\end{theorem}

If higher moments can be bounded, the same property can be established by checking the following criterion.

\begin{cor}\label{cor:KC}
	Let $Y_{\vt}$ and $a_{\vt}$ as in Theorem \ref{thm:KC}.
	Suppose that for each $p>1$ and each $\delta>0$, there exists an $\epsilon>0$ and a constant $C(p,\epsilon)$ such that for $\vs,\vt \in[0,T]^d$, $\|\vs-\vt\|\leq \epsilon$, it holds 
	\begin{align}
		\E  \left| \frac{Y_{\vs}-Y_{\vt}}{\|\vs-\vt\|^{\underline{a}_{[\vs,\vt]}-\delta}} \right|^p \leq C(p,\epsilon, \delta). \label{eqn:moment-condition-allp}
	\end{align}
	Then there exists a modification $\tilde{Y}_{\vt}$ of $Y_{\vt}$ such that \eqref{eqn:KC-result} holds.
\end{cor}

\begin{remark}
	
	The measurability of $\underline{a}_B = \inf_{\vt\in B} a_\vt$ is nontrivial because the infimum might be uncountable. 
	A sufficient condition to ensure measurability is separability of the process $a_\vt$ (see \cite[2.3.3]{vdvWellner1996}): 
	There exists a countable dense subset $\mathcal{G}\subset [0,T]^d$ and an event $\mathcal{N}$ with $P(\mathcal{N})=0$ such that for all $\omega\in\mathcal{N}^c$, and for all $\vt\in [0,T]^d$, there exists a sequence $\vt_n\in\mathcal{G}$ such that $a_{\vt_n}\to a_{\vt}$ as $n\to\infty$.
	The sequence may depend on $\omega$.
	For example, if $d=1$, all cadlag processes are separable. 
	We show in Proposition \ref{prop:measurability} in the appendix that if $a_\vt$ is separable and lower semicontinuous, then $\underline{a}_B$ is measurable for all open and closed sets $B\subset[0,T]^d$, such that Theorem \ref{thm:KC} is applicable.
	As a special example, both conditions are satisfied if $a_\vt$ is continuous.

	If $a_\vt$ is discontinuous, the lower semicontinuity basically requires that at each point of discontinuity, the random field attains the lower value.
	In particular, if condition \eqref{eqn:moment-condition} holds for some random field $a_\vt$, it also holds for the canonical lower semicontinuous variant $a^*_\vt = \lim_{\epsilon\to 0} \inf_{\|\vs-\vt\|<\epsilon} a_\vs$.
	 Moreover, if $a_\vt$ is separable, then $a^*_\vt$ is separable as well, see Proposition \ref{prop:seplsc} in the appendix. 
	Hence, we may conclude that \eqref{eqn:KC-result} holds with $a_\vt$ replaced by $a^*_\vt$.
\end{remark}

The result \eqref{eqn:KC-result} yields that, almost surely, the paths of $\tilde{Y}$ are almost $\underline{a}_B$ Hölder continuous on any closed set $B\subset[0,T]^d$.
Note that this holds pathwise, as $\underline{a}_B$ is possibly random.
Furthermore, this local smoothness holds simultaneously for all closed sets $B$. 
This allows us to study the pointwise Hölder exponent $\alpha_{\vt}(\tilde{Y})$ given by \eqref{eqn:def-pointwise-exponent}, which can be defined in the same way for random fields.
In particular, \eqref{eqn:KC-result} establishes that almost surely, 
\begin{align*}
	\alpha_{\vt}(\tilde{Y}) \geq \lim_{\epsilon\to 0}\inf_{\|\vr-\vt\|<\epsilon} a_{\vr},\quad \forall \vt \in [0,T]^d.
\end{align*}
If the field $a_{\vt}$ itself is continuous or lower semicontinuous, this simplifies to 
\begin{align}
	\alpha_{\vt}(\tilde{Y}) \geq a_{\vt},\quad \forall \vt\in [0,T]^d. \label{eqn:alpha-lower}
\end{align}
Thus, $a_{\vt}$ is a lower bound on the pointwise Hölder exponent. 
We highlight that this lower bound holds with probability $1$ uniformly in $\vt$.
To obtain an equality in \eqref{eqn:alpha-lower}, further properties of the process $Y$ resp.\ $\tilde{Y}$ are required.
For the special case of an Itô-mBm, we establish the equality in Theorem \ref{thm:holder-sharp} below.

\section{A multifractional Gaussian process}\label{sec:mbm}

The aim of this section is to apply the results of the previous section to the Itô-mBm $K^H_t$ employing a random but adapted functional Hurst exponent and to study some properties of the process. 
But before proceeding, let us recall some notions. 
A stationary process has long memory, if its correlations are not integrable. If its correlation function, $ c(\cdot) $, exhibits a power-law behavior at infinity, $ c(\tau) \sim |\tau|^{-\gamma} $, for some $ \gamma \in (0,1) $, the process has long memory with Hurst exponent $ H = 1-\gamma/2 $. Thus, for $ 1/2 < H < 1 $ the fBm has long memory and Hurst parameter $ H $. 
The local properties of the correlations near $0$, however, determine the regularity of the paths and the fractal dimension.
If $ c(\tau) = 1 - b |\tau|^\alpha $, $ \tau \to 0 $, for some $ 0 < \alpha \le 2 $, then, in the sense of \cite{KentWood1997}, the process is locally self-similar of order $ \alpha/2 $, also called similarity index, and the paths have fractal dimension $ D = d+1-\alpha/2 $. 
For Gaussian processes with stationary increments, such as fBm, the local behavior of the incremental variance function (twice the variogram) $V^2(t) = E(X_{s+t}-X_s)^2$, given by $ 2[c(0)-c(t)] $ for a stationary process, characterizes the Hölder exponent and fractal dimension via its expansion in terms of powers of $ |t|^{\alpha} $, see \cite{Orey1970} and \cite{Adler1981}. 
Since for fBm $ B(\cdot,H) $ it holds $ B(x,H) - B(0,H) \stackrel{d}{=}  a^{-H}( B(ax, H) - B(0,H)) $ for all scales $ a> 0 $, local and global properties are linked and the parameter $H$ appears as the Hurst coefficient describing long-memory as well as the self-similiarity index. 
In general, the Hurst parameter as defined above and the local self-similiarity are not necessarily related \cite{GneitingSchlather2004}.
One of the motivations to study multifractional Brownian motions is to have processes which inherit the local behavior of fBm, decoupled from possible long memory effects.
In the sequel, we consider a class of multifractional processes where local and global properties may be decoupled by a suitable choice of an integral kernel, and which contains as a special case the Itô-mBm $K^H_t$ with a time-varying and random exponent $H_s$.
Despite the above discussion on the Hurst coefficient and the self-similarity index, we follow part of the literature and call $ H_s $ (functional) Hurst parameter.

Let $(\Omega, (\F_t)_{t\in\R}, \F, P)$ be a filtered probability space, and let $(W_s)_{s\in\R}$ be a standard Brownian motion w.r.t.\ $\F_t$. 
Furthermore, for each $t\geq 0$, let $g_s(t)$ be $\F_s$-adapted such that $g_s(t)=0$ for $s>t$.
We study the process $X_t$ given by 
\begin{align}
	X_t = \int_{-\infty}^{t} g_s(t)\, dW_s, \label{eqn:def-X}
\end{align}
where the kernel $g_s(t)$ satisfies $\int_{-\infty}^{t} |g_s(t)|^2\, ds < \infty$ for all $t$, such that the stochastic integral is well-defined, see \cite[Section 3.2.D]{Karatzas1998}.  
Furthermore, we assume that the integrand is of the following form.

\begin{condition*}[A]\label{cond:g_s(t)}
	The function $t\mapsto g_s(t)$ is differentiable in $t>s$ for all $s$. 
	There exist $\F_s$-adapted processes $H_s$, $L_s$, and $R_s$, such that $H_s\in(0,1)$, $R_s>\frac{1}{2}$, and it holds for all $t\geq 0$,
	\begin{align}
		\left| g_s(t)\right| &\leq   L_s |t-s|^{H_s-\frac{1}{2}},&&\qquad s\in(t-1,t), \label{eqn:cond-B-3} \\
		\left|\partial_t g_s(t)\right| 		&\leq L_s|t-s|^{H_s-\frac{3}{2}} ,&&\qquad s\in(t-1,t). \label{eqn:cond-B-1} \\
		\left| \partial_t g_s(t) \right| 	&\leq L_s|t-s|^{-R_s}, &&\qquad s\in(-\infty,t-1]. \label{eqn:cond-B-2}
	\end{align}
\end{condition*}

\begin{figure}[tb]
	\includegraphics[width=0.48\textwidth]{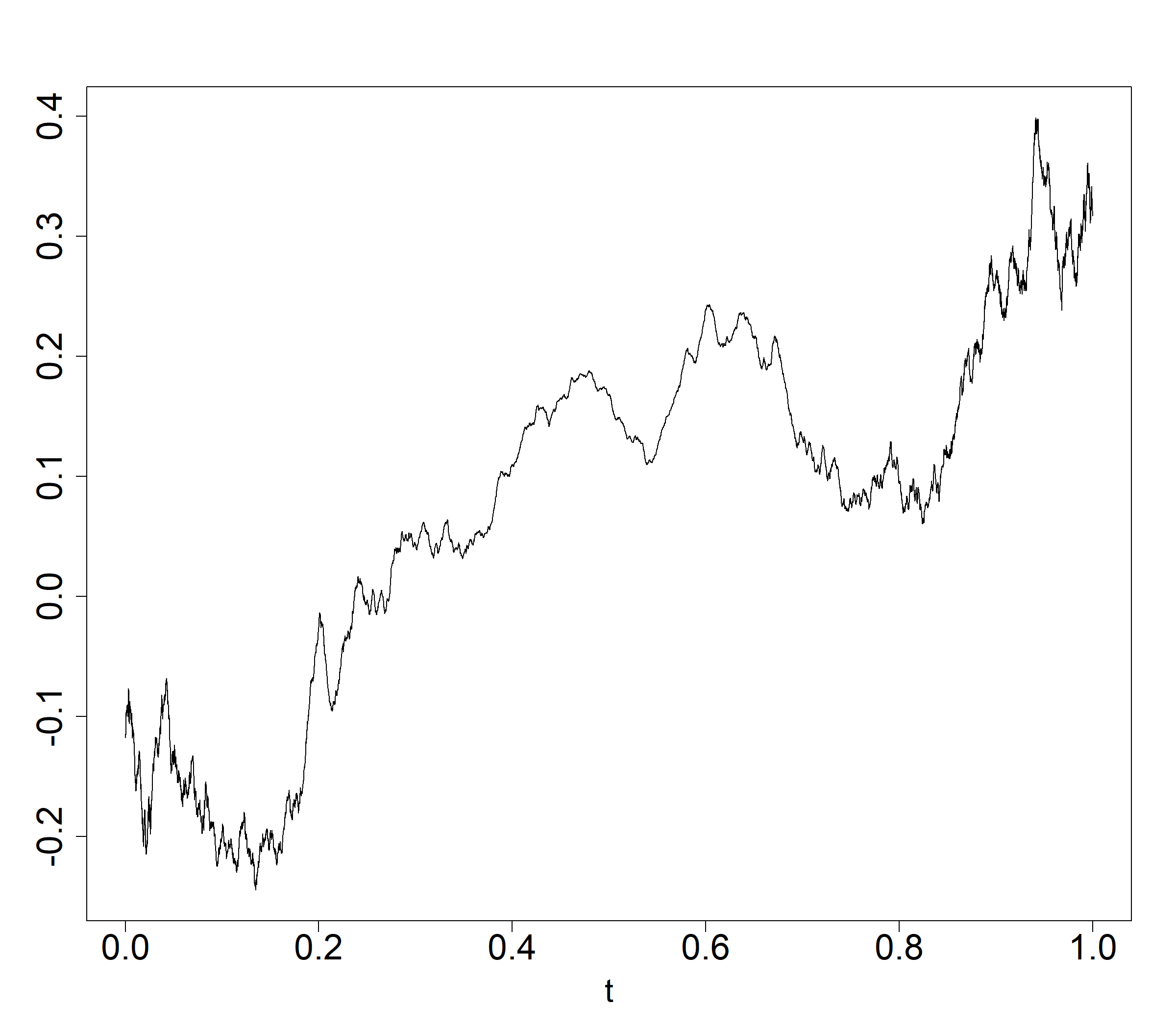}\includegraphics[width=0.48\textwidth]{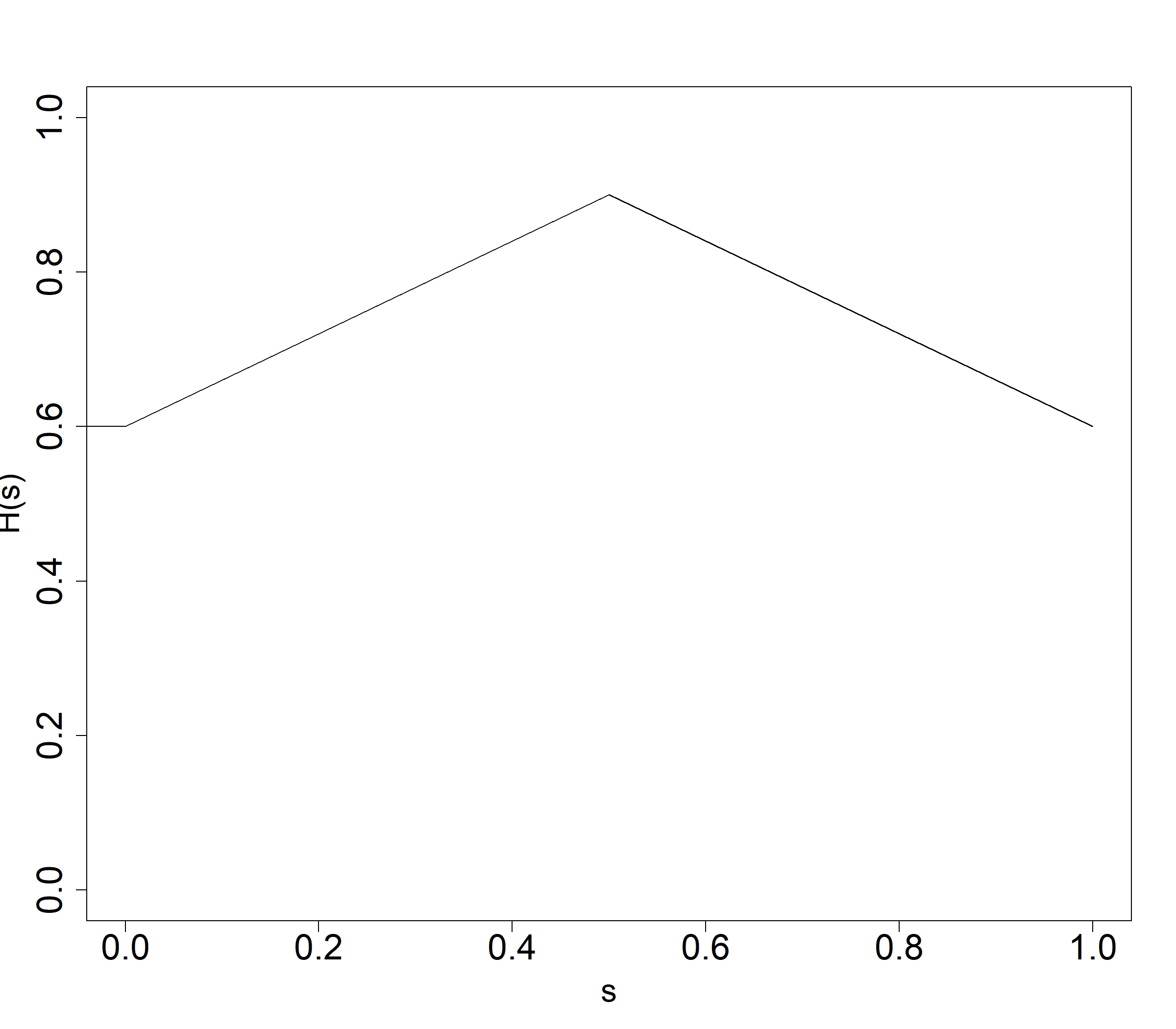}
	\caption{Left: sample path of the multifractional Matérn process with kernel \eqref{eqn:matern} and parameter $\lambda=4$. Right: corresponding path of $H_s$.}
	\label{fig:matern}
\end{figure}

We will restrict ourselves to stochastic processes indexed by a univariate time parameter $t$.
While the primary example to satisfy Condition \nameref{cond:g_s(t)} is the Itô-mBm \eqref{eqn:mbm-ito}, our condition allows for more general processes.
For example, we may introduce a logarithmic factor and consider
\begin{align*}
	g_s(t) =  \left[ (t-s)_+\log(t-s)_+\right]^{H_s-\frac{1}{2}} - \left[(-s)_+ \log(-s)_+\right]^{H_s-\frac{1}{2}}.
\end{align*}
Moreover, short term and long term behavior may be decoupled, e.g.\ by defining 
\begin{align*}
	g_s(t) = (t-s)_+^{H_s-\frac{1}{2}} \phi(|s-t|)
\end{align*}
for a bounded function $\phi$ such that $\phi(x)=1$ near zero, and $\phi(x)\to 0$ rapidly as $x \to \infty$.
If $\phi$ is compactly supported, the process $X_t$ will have a finite memory.
As a further example, let us consider a generalization of the well known Matérn process. The Matérn process, see \cite{Lilly2017}, is given by
\begin{align*}
	g^M_s(t) = \frac{1}{\Gamma(H+1/2)} (t-s)_+^{H-\frac{1}{2}} e^{- \lambda (t-s)},
\end{align*}
for a parameter $ \lambda > 0 $ which has the interpretation of a damping factor, as the spectrum is flat and white-noise like for frequencies $ \bar{\omega} \ll \lambda $ and has a power-law decay for high frequencies $ \bar{\omega} \gg \lambda $. For $ |\tau| \gg 1/\lambda $ the autocovariance function, $ \gamma^M(\tau ) $, behaves like $ | \lambda \tau |^{H-1/2} e^{-\lambda|\tau|} $, such that the process has short memory. For small scales $ |\tau | \ll 1/\lambda $ one has $ \gamma^M( \tau ) \approx \sigma^2 - M_H |\tau|^{2H} $, where $ \sigma^2 $ denotes the variance and $M_H$ is a constant. This means, the Matérn process is (only) locally self-similar, see \cite{KentWood1997}, whereas fBm is self-similar at all scales. Let us introduce the following generalization given by
\begin{align}
	g_s(t) = (t-s)_+^{H_s-\frac{1}{2}} e^{- \lambda (t-s)}. \label{eqn:matern}
\end{align}
for an $\mathcal{F}_s$-adapted process $H_s\in(0,1)$.
We will call the process corresponding to the kernel \eqref{eqn:matern} the multifractional Matérn process.
A sample path is depicted in Figure \ref{fig:matern}, which has been obtained by discretizing the integrand on the interval $[-10,1]$ with mesh size $10^{-5}$. 
As we will see below, the local regularity of $X_t$ is determined by $H_t$.
The generalized Matérn process and the above discussion illustrates that the class of processes satisfying Condition (A) is large and covers short memory as well as long memory processes with locally varying and random self-similarity index $H_t$.

To complement the assumptions on $g_s(t)$, we require the following conditions.

\begin{condition*}[B]\label{cond:boundedness}
	There exist real numbers $\underline{H}$, $\overline{H}$, $\overline{L}$, and $\underline{R}$, such that for all $t\geq 0$
	\begin{align}
		\begin{split}
			0<\underline{H}\leq H_{t} \leq \overline{H} <1, \\
			|L_{t}|\leq \overline{L}, \quad\text{and}\quad |R_{t}|\geq \underline{R}>\frac{1}{2}.
		\end{split} \label{eqn:global-boundedness}	
	\end{align}
\end{condition*}

\begin{condition*}[C]\label{cond:continuity}
	There exists a continuous, increasing function $\omega(h)$ with $\omega(0)=0$, such that almost surely, it holds for all $t\geq 0$ and $h>0$ that  
	\begin{align}
		\begin{split}
			|H_{t+h} - H_{t}| \leq \omega(h),
		\end{split}\label{eqn:smooth-H} 
	\end{align}
\end{condition*}

Condition \nameref{cond:boundedness} imposes a convenient type of boundedness, and Condition \nameref{cond:continuity} requires that $H_t$ is uniformly continuous with a deterministic bound on the modulus of continuity.
The assumptions on $\omega$ are minimal. 
Condition \nameref{cond:continuity} is satisfied, for instance, if $H_t$ is a deterministic continuous function which is not even Hölder continuous.

While conditions \nameref{cond:boundedness} and \nameref{cond:continuity} might seem restrictive for random integrands $g_s(t)$, they can be relaxed by a localization procedure \cite[Sec.\ 4.4.1]{jacod2011discretization}. 
For a sequence of stopping times $\tau_n\to\infty$, the stopped process $X_{t\wedge\tau_n}$ may be represented as \eqref{eqn:def-X} with integral kernel $g^n_s(t) = g_s(t\wedge \tau_n)$.
Nevertheless, pathwise regularity of the stopped processes $X_{t\wedge\tau_n}$  on the finite interval $(0,T)$ readily transfers to $X_t$ because $\tau_n\geq T$ eventually. 
Hence, instead of conditions \nameref{cond:boundedness} and \nameref{cond:continuity}, we could alternatively impose the following localized assumptions.

\begin{condition*}[L-BC]\label{cond:continuity-local}
	There exists a sequence of stopping times $\tau_n\to\infty$, real numbers $\underline{H}^n$, $\overline{H}^n$, $\overline{L}^n$, $\underline{R}^n$, and a sequence of continuous, increasing functions $\omega_n(h)$ with $\omega_n(0)=0$, such that \begin{align}
		\begin{split}
			|H_{(t+h)\wedge \tau_n} - H_{t\wedge\tau_n}| \leq \omega_n(h),
		\end{split}\label{eqn:smooth-H-local} \\
		\begin{split}
			0<\underline{H}^n\leq H_{t\wedge\tau_n} \leq \overline{H}^n <1, \\
			|L_{t\wedge\tau_n}|\leq \overline{L}^n, \text{and}\; |R_{t\wedge\tau_n}|\geq \underline{R}^n>\frac{1}{2}.
		\end{split} \label{eqn:local-boundedness}	
	\end{align}
\end{condition*}

The local boundedness assumption \eqref{eqn:local-boundedness} is rather mild, and it is satisfied whenever the processes are continuous, or have bounded jumps.
Furthermore, while \eqref{eqn:smooth-H-local} requires $H_t$ to be continuous, we do not impose any Hölder regularity as the restrictions on $\omega_n$ are minimal.
In fact, a suitable sequence of stopping times can be constructed if only $H_t$ is $\eta_T$-Hölder continuous on $[0,T]$ for any $T$. 
The value $\eta_T>0$ might even be random.
To be precise, for two sequences $\gamma_n\downarrow 0$, $M_n\to \infty$, we could define
\begin{align*}
	\tau_n = \inf \left\{ t \geq 0: \sup_{r,s\in(0,t)} \frac{|H_r-H_s|}{|r-s|^{\gamma_n}} > M_n \right\}.
\end{align*}
Clearly, $\tau_n$ are stopping times such that condition \nameref{cond:continuity-local} holds, and $\tau_n\to\infty$ holds if $H_t$ is Hölder continuous as in the current situation. 
For instance, the continuity condition \eqref{eqn:smooth-H-local} is satisfied by a continuous Itô process, or fractional Brownian motion with any Hurst exponent. 
In this case, the local boundedness of $H_t$ may be achieved by a repelling drift term.

To simplify the proofs in the following, we will work without loss of generality with the stronger assumptions \nameref{cond:boundedness} and \nameref{cond:continuity}.
Using the previous assumptions, we can establish the following moment inequalities, which are a prerequisite for the application of Theorem \ref{thm:KC}.

\begin{lemma}\label{lem:inc-moment}
	Let conditions \nameref{cond:g_s(t)}, \nameref{cond:boundedness}, and \nameref{cond:continuity} hold, and let $\epsilon\in(0,1)$. 
	Then for all $h\in(0,\epsilon)$, $\delta>0$, any $t\geq 0$, and any $p>1$,
	\begin{align}
		\E \left| \frac{X_{t+h} - X_t}{h^{H_t}} \right|^p \leq C(\underline{H}, \overline{H}, \overline{L}, \underline{R},\epsilon,p, \delta)\,h^{-2p\,\omega(2\epsilon)-p\delta}. \label{eqn:inc-moment-res-1}
	\end{align}
	Alternatively, allowing $H_t$ to be discontinuous such that conditions \nameref{cond:g_s(t)} and \nameref{cond:boundedness} hold,
	\begin{align}
		\E \left| \frac{X_{t+h} - X_t}{h^{\underline{H}}} \right|^p \leq C(\underline{H}, \overline{H}, \overline{L}, \underline{R},\epsilon,p, \delta) h^{-p\delta}. \label{eqn:inc-moment-res-2}
	\end{align}
	Finally, suppose $H_t$ is deterministic and potentially discontinuous, such that conditions \nameref{cond:g_s(t)} and \nameref{cond:boundedness} hold. 
	Denote $\underline{H}^\epsilon_t=\inf_{s\in[t-\epsilon,t+\epsilon]} H_s$. 
	Then 
	\begin{align}
		\E \left| \frac{X_{t+h} - X_t}{h^{\underline{H}^\epsilon_t}} \right|^p \leq C(\underline{H}, \overline{H}, \overline{L}, \underline{R},\epsilon,p,\delta) h^{-p\delta}. \label{eqn:inc-moment-res-3}
	\end{align}
\end{lemma}

These moment bounds can be used to derive the local regularity of the paths of $X_t$ via the continuity criterion of Theorem \ref{thm:KC}.

\begin{theorem}\label{thm:regularity-X}
	Let conditions \nameref{cond:g_s(t)} and \nameref{cond:boundedness} hold. 
	If $H_t$ satisfies condition \nameref{cond:continuity}, then there exists a version of $X$, also denoted as $X$, such that
	\begin{align}
		P\left( \alpha_t(X)\geq H_t \; \forall t\in(0,T) \right) = 1. \label{eqn:smoothness-local}
	\end{align}
	Alternatively, if $H_t$ is potentially discontinuous, then
	\begin{align}
		P\left( \alpha_t(X)\geq \underline{H} \; \forall t\in(0,T) \right) = 1. \label{eqn:smoothness-lower}
	\end{align}
	If $H_t$ is potentially discontinuous but deterministic, then
	\begin{align}
		P\left( \alpha_t(X)\geq H_t^* \; \forall t\in(0,T)  \right) = 1, \label{eqn:smoothness-deterministic} 
	\end{align}
	where $H_t^* = \lim_{\epsilon\to 0} \inf_{|r-t|<\epsilon} H_r$ is the lower semicontinuous variant of $H_t$.
\end{theorem}

This regularity result can be compared to the properties derived by \cite{Ayache2018}.
Therein, it is shown that the Hölder exponent $\alpha(X, [s,t])$ on any interval $[s,t]$ satisfies $\alpha(X, [s,t]) \geq \inf_{r\in[s,t]} H_r$ almost surely, where
\begin{align*}
	\alpha(X,[s,t]) = \sup\left\{ \gamma\in[0,1]: \sup_{r_1,r_2\in [s,t]} \frac{|X_{r_1}-X_{r_2}|}{|r_1-r_2|^\gamma} <\infty \right\} .
\end{align*}
The same property is established in the proof of Theorem \ref{thm:regularity-X} above, see \eqref{eqn:Ito-mBm-uniform-Holder}, which also establishes the inequality uniformly in $s,t$.
Furthermore, \cite{Ayache2018} requires $H_t>\frac{1}{2}$ and that $H_t$ is $\eta$-Hölder continuous for some $\eta>\frac{1}{2}$.
In contrast, we do not restrict the range of $H_t$, and our continuity assumption on $H_t$ is much weaker.
As a result, the regularity of $X_t$ is mostly decoupled from the regularity of $H_t$.

\begin{figure}[tb]
	\centering
	\includegraphics[width=0.48\textwidth]{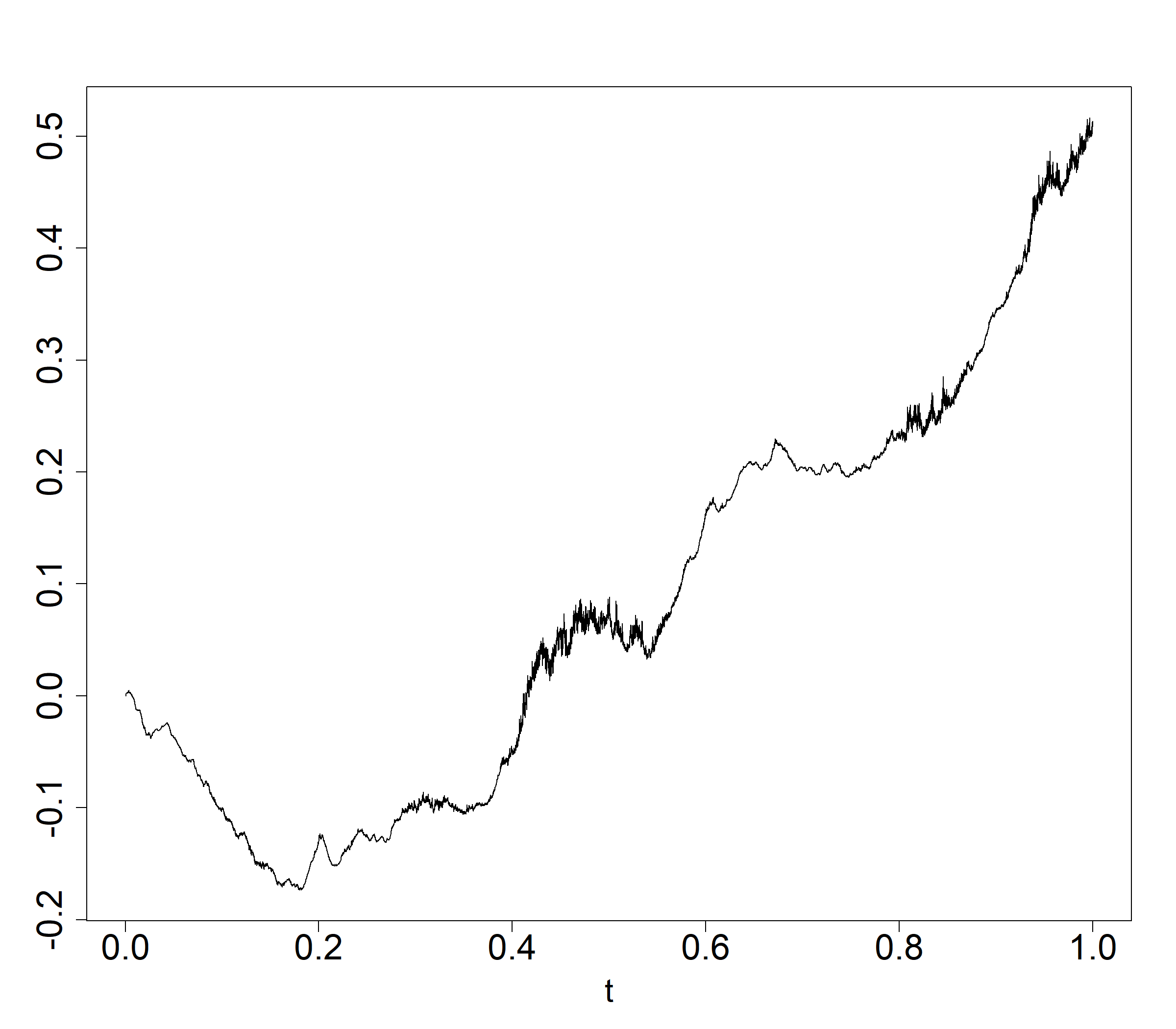}
	\includegraphics[width=0.48\textwidth]{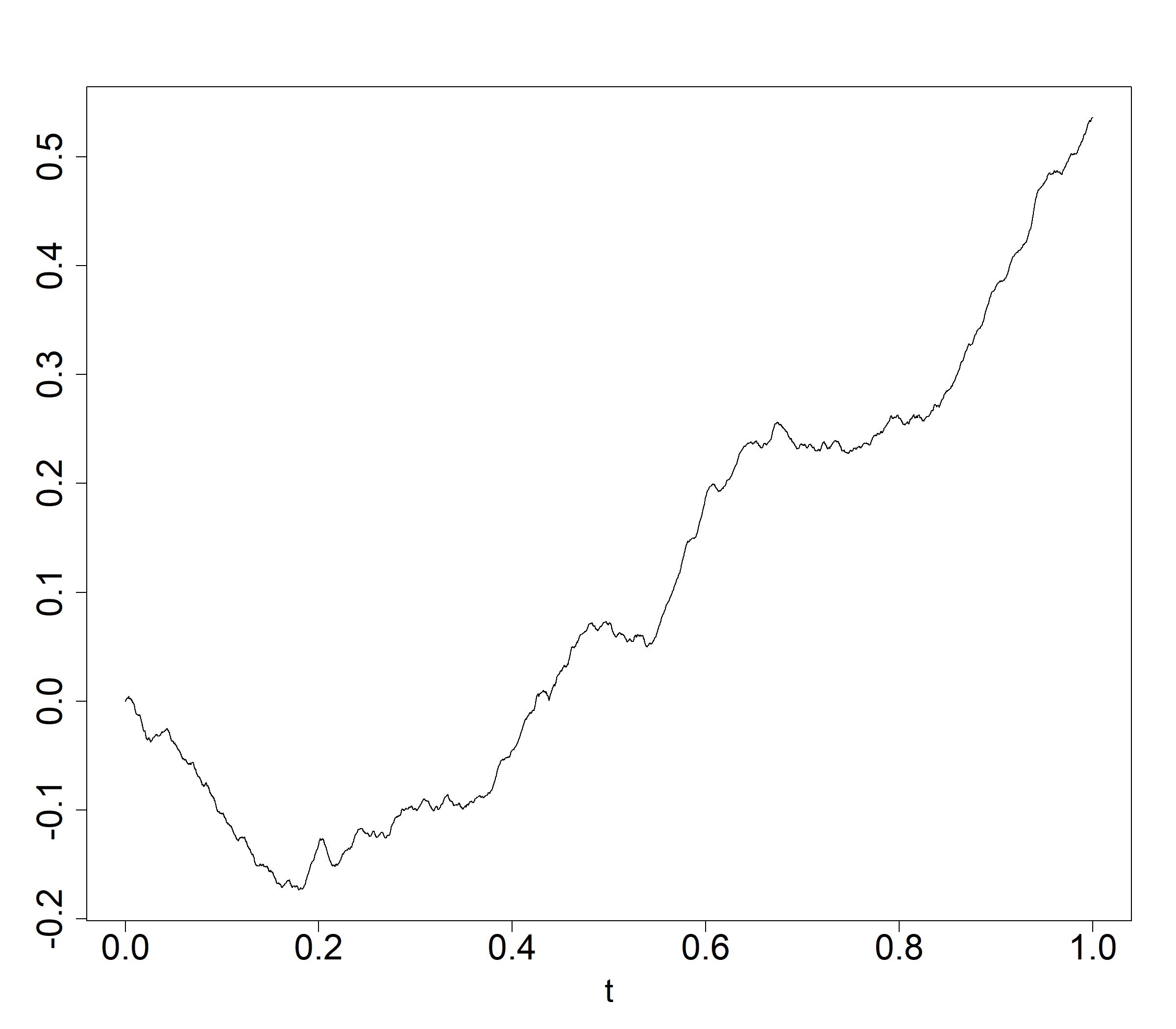}
	\caption{Sample paths of multifractional Brownian motion (left) and Itô-mBm (right) for $\sigma_t=1$ and a rough local Hurst parameter $H_t\in(0.85,0.95)$ constructed as $H_t = 0.9 + 0.05\tanh(B^{0.2}_t)$.}
	\label{fig:mbm-rough}
\end{figure}

This finding is illustrated by Figure \ref{fig:mbm-rough}, where we depict sample paths of the classical mBm $B_t^H$ and the Itô-mBm $X_t=K_t^H$, for the same functional Hurst parameter $H_t$. 
Here, the Hurst parameter is constructed as $H_t = 0.9 + 0.05\tanh(\tilde{B}_t^{0.2})$ for a fractional Brownian motion $\tilde{B}_t^{0.2}$ with Hurst parameter $0.2$, independent of $W_t$. 
In particular, $H_t\in(0.85, 0.95)$, and $H_t$ is almost $0.2$-Hölder continuous, which is rather rough.
The plot demonstrates that the process $X_t$ is almost $0.85$-Hölder continuous, irrespective of the roughness of $H_t$. 
In contrast, $B_t^H$ is as rough as $H_t$. 
To simulate the sample paths, the integrand of the stochastic integral has been discretized on the interval $[-10,1]$ with mesh size $10^{-5}$. 
The code used for the simulations is available as a digital supplement. 

It is natural to ask whether equality holds in \eqref{eqn:smoothness-local} and \eqref{eqn:smoothness-deterministic}.
In order to obtain upper bounds on the Hölder smoothness of $X_t$, we study the local behavior of increments $X_{t+h}-X_t$ as $h$ tends to zero.
For the commonly studied multifractional Brownian motion $B_t^H=B(t,H_t)$, it is known that \begin{align*}
	h^{-H_t}\left[ B(t+hr, H_{t+hr})-B(t,H_t) \right] \wconv B(r, H_t),
\end{align*}
whenever $H_t$ is a smooth deterministic function.
This can be seen, e.g., by using that $H\mapsto B(t,H)$ is $C^\infty$ \cite{Ayache2013}.
Hence, locally around $t$, the mBm behaves like a fractional Brownian motion with Hurst parameter $H_t$.
Intuitively, this bounds its pointwise Hölder exponent from above by $H_t$. 
For our process $X_t$, the same scaling behavior holds, as established in the sequel.
To this end, we need to impose further restrictions on the integral kernel $g_s(t)$, which could also be localized by a sequence of stopping times. 

\begin{condition*}[A*]\label{cond:B2}
	There exists a continuous process $\sigma_t$ satisfying $|\sigma_t|\leq L_t$ for the same $L_t$ as above, and a $\rho>0$, such that for all $t\geq 0$
	\begin{align*}
		\left| g_s(t) - \sigma_s (t-s)^{H_s-\frac{1}{2}} \right| & \leq L_s |t-s|^{H_s-\frac{1}{2}+\rho}, \qquad s\in(t-1,t), \\
		\left| \partial_t g_s(t) - \partial_t \sigma_s (t-s)^{H_s-\frac{1}{2}} \right| &\leq L_s |t-s|^{H_s-\frac{3}{2}+\rho}, \qquad s\in(t-1,t).
	\end{align*}
\end{condition*}

Note that condition \nameref{cond:B2} amounts to imposing assumption \nameref{cond:g_s(t)} for the difference kernel $g_s(t) - \sigma_s |t-s|^{H_s-\frac{1}{2}}$, with slightly larger exponents.
The small $\rho>0$ is required to ensure that the approximation error is asymptotically negligible.

The following limit theorem establishes weak convergence of the locally rescaled process to a fractional Brownian motion.
This convergence is functional, i.e.\ it holds in the space $C(I)$ of continuous functions on an interval $I$.
Furthermore, we show that the convergence holds $\F_{t-}$-stably in distribution \cite{Podolskij2010}.
This property is not used in the sequel, but might be interesting for future statistical applications.

\begin{theorem}[Rescaling limit]\label{thm:rescaling}
	Suppose that conditions \nameref{cond:g_s(t)}, \nameref{cond:B2}, \nameref{cond:boundedness}, and \nameref{cond:continuity} hold.
	Assume furthermore that the modulus of continuity satisfies $\omega(\epsilon) \log(\epsilon) \to 0$ as $\epsilon\to 0$.
	Let $t>0$ be fixed, $a<0<b$ and $I=(a,b)$, and $Z$ be any bounded random variable which is $\F_{t-}$-measurable. 
	Then
	\begin{align*}
		\begin{pmatrix}
			h^{-H_t} \left(X_{t+hr}-X_t\right)_{r\in I}
			\\ Z
		\end{pmatrix}
		\wconv 
		\begin{pmatrix}
			\left(\displaystyle\int_{-\infty}^{r} \sigma_t[(r-z)_+^{H_t-\frac{1}{2}} - (-z)_+^{H_t-\frac{1}{2}}]\, d\tilde{W}_z \right)_{r\in I}
			\\ Z
		\end{pmatrix}
	\end{align*}
	as $h\to 0$, where $\tilde{W}_z$ is a Brownian motion independent of $\F_{t-}$.
	Convergence holds in $C(I)\times \R$. 
\end{theorem}

In view of the previous rescaling limit, the Itô-mBm may be considered as a non-stationary generalization of fBm.
Note that the condition $\omega(\epsilon) = o(\log \epsilon)$ is also required by \cite[Thm.\ 3.2]{Surgailis2008}, where a rescaling limit for a related multifractional process is derived.
The processes studied therein are defined in terms of an integral representation similar to \eqref{eqn:mbm-ito}, but they are not contained in our framework \eqref{eqn:def-X}.

As the weak convergence in Theorem \ref{thm:rescaling} is functional, we are able to derive an upper bound on the pointwise Hölder exponent at time point $t$ of the paths of $X$.

\begin{theorem}\label{thm:holder-sharp}
	If the conditions of Theorem \ref{thm:rescaling} hold, then for each $t>0$ and any continuous version of $X_t$, we have 
	\begin{align*}
		P\left(\alpha_t(X) = H_t\right)=1.
	\end{align*}
	Moreover, we may describe the Hölder exponent uniformly on all intervals $(s,t)$ by
	\begin{align*}
		P\left( \forall s,t\in[0,T), s<t \,:\, \inf_{r\in[s,t]} H_{r} = \sup \left\{ \alpha\,:\, \lim_{\epsilon\to 0} \sup_{r,v\in (s,t), |r-v|<\epsilon} \frac{|X_r-X_v|}{|r-v|^\alpha} =0 \right\} \right) =1.
	\end{align*}
\end{theorem}

Let us briefly discuss some implications and applications of the above results for a Hölder-continuous process $X_t$ of regularity $\alpha = \underline{H}$. 
As $\alpha$-Hölder continuous functions have bounded $p$-variation, $p = 1/\alpha$, for any function $h$ with bounded $q$-variation, $1/q+\alpha >1$, the integrals $\int_0^t h_{s-} \, d X_s$ and $ \int_0^t X_{s-} \, d h_s $ exist a.s.\ as pathwise refinement of Riemann Stieltjes integrals (also called Young-type integral), and the associated integral operator is continuous, since 
\begin{align*}
	\left| \int_0^t h_s \, d X_s - h_0 (X_t-X_0) \right| \le C \| X \|_{[1/\alpha],[0,t]} \| h \|_{[q],[0,t]}
\end{align*}
for some constant $ C = C(\alpha,q) $, see \cite{Dudley1992} and \cite{Norvaisa2006}. 
Here $ \| \cdot \|_{[r],[0,t]} $ is the $r$-variation norm w.r.t. the interval $ [0,t] $. 
The integrand $h_s$ can also be a càdlàg process with paths of finite $q$-variation.  Morever, the Young inequalities
\begin{align*}
	\left| \int_s^t h_{s-} \, d X_s - h_{s-}( X_t - X_s ) \right| \le C \| h_{-} \|_{[q],[s,t]} \| X \|_{[p],[s,t]}, \\
	\left| \int_s^t h_{s-} \, d X_s - h_{s}( X_t - X_s ) \right| \le C \| h \|_{[q],[s,t]} \| X \|_{[p],[s,t]}
\end{align*}
hold. 

Fractional Brownian motion and multifractional Brownian motion are often studied in terms of their covariance function.
The covariance function of fBm as in \eqref{eqn:def-fbm} is given by \cite[Thm.\ 1.3.1]{Mishura2008}
\begin{align}
	\E\left[B(t,H)\, B(s,H)\right] &= \frac{ A(H)}{2} \left( |t|^{2H}+|s|^{2H} - |t-s|^{2H} \right), \label{eqn:cov-fbm}\\
	A(H) &= \frac{\Gamma(H+\tfrac{1}{2})^2}{2 H \sin(\pi H) \Gamma(2H)},\quad H\in(0,1). \nonumber
\end{align}

The process $X_t$ given by \eqref{eqn:def-X} is not fully specified, hence the complete covariance function can not be determined in closed form.
However, we may obtain an explicit limit expression for the covariance when rescaling the increments of the process.
In fact, the proof of Theorem \ref{thm:rescaling} reveals that in addition to the weak convergence, we also obtain convergence of the local covariance function.

\begin{cor}\label{cor:covariance-local}
	Under the conditions of Theorem \ref{thm:rescaling}, for any fixed $t>0$, the covariance function $\gamma_{t,h}$ of the rescaled process satisfies,
	\begin{align}
		\gamma_{t,h}(r,v)&= \E\left[ h^{-H_t} \left( X_{t+hr}-X_t \right) h^{-H_t} \left( X_{t+hv}-X_t \right) \right] \nonumber \\
		&\to   \E \left[\sigma_t^2\, \frac{A(H_t)}{2}\, \left( |r|^{2H_t} + |v|^{2H_t}-|r-v|^{2H_t} \right) \right], \label{eqn:local-cov}
	\end{align}
	for any $r,v\in\R$, as $h\to 0$.
\end{cor}

Equation \eqref{eqn:local-cov} yields that conditionally on $H_t$ the process $X_t$ is locally self-similar in the sense of \cite[Section 2]{Cavanaugh2003}. 
If $H_t$ is deterministic, \eqref{eqn:local-cov} is precisely the covariance function of fractional Brownian motion with Hurst parameter $H=H_t$.

In order to study the long-term covariances of $X_t$, we need to specify the integrand $g_s(t)$ in \eqref{eqn:def-X} more explicitly.
To this end, we investigate the Itô-mBm $K_t^H$ given by \eqref{eqn:mbm-ito}. 
Itô's isometry yields that, for $s,t\in\R$,
\begin{align}
	&\qquad \E \left( K_t^H, K_s^H \right) \nonumber\\
	&= \E \int_{-\infty}^{\infty} \sigma_r^2 \left[(t-r)_+^{H_r-\frac{1}{2}} - (-r)_+^{H_r-\frac{1}{2}}\right] \left[(s-r)_+^{H_r-\frac{1}{2}} - (-r)_+^{H_r-\frac{1}{2}}\right] \, dr. \label{eqn:cov-Itombm}
\end{align}
This expression is inconvenient as it depends on the whole path of $H_t$.
In contrast, for the classical mBm, it can be readily checked that the covariance $\E[B_t^H\, B_s^H]$ depends on the function $H$ only via the values $H_s$ and $H_t$.
Definition \eqref{eqn:mbm-field} with deterministic Hurst function corresponds to the process $Y_{(1,0)}$ as considered by \cite{Stoev2006}, hence Theorem 4.1 therein demonstrates that 
\begin{align}
	\begin{split}
		\E\left[ B_t^H\, B_s^H \right]
		&=D(H_t, H_s) \left[ |t|^{2H_{t,s}}\cos( \pi(\tilde{H}_{t,s} - \sgn(t)H_{t,s})  ) \right. \\
		&\qquad \left.+ |s|^{2H_{t,s}}\cos( \pi(\tilde{H}_{t,s} + \sgn(s)H_{t,s})  ) \right. \\
		&\qquad \left.- |t-s|^{2H_{t,s}}\cos( \pi(\tilde{H}_{t,s} - \sgn(t-s)H_{t,s})  ) \right], 
	\end{split} \label{eqn:cov-mbm}\\
	D(H_t, H_s) 
	&= \frac{\Gamma(H_t+\tfrac{1}{2})\, \Gamma(H_s+\tfrac{1}{2})\, \Gamma(2-2H_{t,s})}{2\pi H_{t,s} (1-2H_{t,s})}, \nonumber\\
	H_{t,s} 
	&= \tfrac{H_s+H_t}{2}, \qquad \Delta \tilde{H}_{t,s} = \tfrac{H_s-H_t}{2}. \nonumber
\end{align}
Note that \eqref{eqn:cov-mbm} reduces to \eqref{eqn:cov-fbm} if $H_t=H_s$.
There exist alternative definitions of mBm for which the covariance function has a slightly simpler form, e.g.\ based on the harmonizable representation of fractional Brownian motion studied in \cite{Ayache2000}.
These alternatives are compared by \cite{Stoev2006}.

Despite its path-dependence, the covariance \eqref{eqn:cov-Itombm} of the Itô-mBm can be simplified in the case where $H_t$ and $\sigma_t$ are stationary processes. 

\begin{proposition}\label{prop:cov-stationary}
	Suppose that $H_t\in(0,1)$ and $\sigma_t>0$ are stationary adapted processes. 
	Then 
	\begin{align*}
		\E \left( K_t^H, K_s^H \right)
		&= \E \left[\sigma_0^2 \, \frac{A(H_0)}{2}\, \left( |t|^{2H_0} +|s|^{2H_0}-|t-s|^{2H_0}  \right)\right], \quad s,t\in\R,
	\end{align*}	 
	provided that the latter value is finite for all $s,t\in\R$.
\end{proposition}

In particular, we may compute the autocovariance function of the toppstationary increment process $(K_{t+1}^H-K_t^H)$ as 
\begin{align*}
	\gamma(\Delta) &= \E \left[ (K_{t+1}^H-K_t^H)(K_{t+\Delta+1}^H-K_{t+\Delta}^H) \right] \\
	&= \E \left[ |\Delta+1|^{2H_0} -2|\Delta|^{2H_0} + |\Delta-1|^{2H_0} \right], \qquad \Delta\in\mathbb{N}.
\end{align*}
This expression could be used to derive the long range dependence of the stationary increments. 

\section*{Acknowledgments}

The authors would like to thank the anonymous reviewer for his careful reading and constructive comments which helped to improve this article considerably.

\section{Proofs}\label{sec:proofs}

\subsection*{Proofs of Section \ref{sec:KC}}

Recall that any lower semicontinuous function $f : \R^d \to \R $,  $\lim_{\epsilon \to 0} \inf_{y: \|y-x\|<\epsilon} f(y) \geq  f(x)$, has a minimimum over any compact set $B$. The proof of Theorem \ref{thm:KC} requires the following preliminary result, which asserts that the minimum over $B$ equals the limiting infimum over its $ \epsilon $-enlargements, as $ \epsilon \to 0 $. 

\begin{lemma}\label{lem:LSC}
	Let $B\subset\R^d$ be a compact set, and let $f:\R^d\to \R$ be lower semicontinuous.
	Then 
	\begin{align*}
		\lim_{\epsilon\to 0}\inf_{x\in B} \inf_{y:\|y-x\|< \epsilon} f(y)=  \inf_{x\in B} f(x).
	\end{align*}
	This limit is monotonically increasing as $\epsilon\downarrow 0$. 
\end{lemma}
\begin{proof}[Proof of Lemma \ref{lem:LSC}]
	Clearly, it suffices to show '$\leq$'.  Let $\eta>0$.
	By the lower semicontinuity, for each $x\in B$, there exists an $\delta_x>0$ such that $f(y)\geq f(x)-\eta$ for all $y\in B_{\delta_x}(x)$. 
	The balls $B_{\delta_x}(x)$ are a covering of the compact set $B$.
	Hence, there exist $x_1,\ldots, x_n\in B$ such that $B\subset \bigcup_{i=1}^n B_{\delta_{x_i}}(x_i)$.
	Furthermore, for $\epsilon=\epsilon(\eta)$ sufficiently small, it holds that $B^\epsilon \subset \bigcup_{i=1}^n B_{\delta_{x_i}}(x_i)$, where $B^\epsilon = \bigcup_{x\in B} B_{\epsilon}(x)$ is the open $\epsilon$-neighborhood of $B$.
	Thus,
	\begin{align*}
		\inf_{x\in B} f(x) \quad\geq\quad \inf_{y\in B^\epsilon} f(y) &\geq \min_{i=1,\ldots, n} \inf_{y\in B_{\delta_{x_i}}(x_i)} f(y) \\
		&\geq \min_{i=1,\ldots, n} f(x_i) - \eta \\
		&\geq \inf_{x\in B} f(x) - \eta.
	\end{align*}
	Since $\eta>0$ is arbitrary, this establishes the claim of the Lemma.
	The monotonicity is evident because $\inf_{x\in B} \inf_{y: \|y-x\|<\epsilon_1} f(y) \leq \inf_{x\in B} \inf_{y: \|y-x\|<\epsilon_2} f(y)$ for $\epsilon_1 \geq \epsilon_2$. 
\end{proof}

Clearly, Lemma \ref{lem:LSC} also holds if we replace the whole space $\R^d$ by some open set, e.g.\ $(0,T)^d$.

\begin{proposition}\label{prop:measurability}
	Let $a_\vt$ be a random field indexed by $\vt\in[0,T]^d$, i.e.\ each $a_\vt$ is a random variable.
	If $a_\vt$ is separable and lower semicontinuous, then for any $B\subset \R^d$ which is open or closed, $\underline{a}_B = \inf_{\vt\in B \cap [0,T]^d} a_\vt$ is a random variable. 
\end{proposition}
\begin{proof}[Proof of Proposition \ref{prop:measurability}]
	Let $\mathcal{G}\subset(0,T)^d$ be the countable dense set such that $a_\vt$ is separable w.r.t.\ $\mathcal{G}$, and let $\mathcal{N}$ be the corresponding null event.
	Separability of $a_\vt$ implies that, on $\mathcal{N}^c$,
	\begin{align*}
		\inf_{\vt\in B \cap [0,T]^d} a_\vt = \inf_{\vt \in B\cap \mathcal{G}} a_\vt,
	\end{align*}
	for any open set $B\subset [0,T]$, which is measurable as a countable infimum.
	If instead $B$ is closed, we use Lemma \ref{lem:LSC} such that 
	\begin{align*}
		\inf_{\vt\in B} a_\vt = \lim_{k\to \infty} \inf_{\vt \in B^{1/k} \cap [0,T]^d} a_\vt,
	\end{align*}
	where $B^\epsilon$ denotes the open $\epsilon$-neighborhood of the set $B$. 
	In particular, $\underline{a}_B$ is measurable.
\end{proof}

\begin{proposition}\label{prop:seplsc}
	Let $a_\vt$ be a random field indexed by $\vt\in[0,T]^d$, which is separable. 
	Then $a^*_\vt = \lim_{\epsilon\to 0} \inf_{\vs\in[0,T]^d,\|\vs-\vt\|<\epsilon} a_\vs$ is as well a separable and lower semicontinuous random field. 
\end{proposition}
\begin{proof}[Proof of Proposition \ref{prop:seplsc}]
	Let $\mathcal{G}\subset[0,T]^d$ be the countable dense set such that $a_\vt$ is separable w.r.t.\ $\mathcal{G}$, and let $\mathcal{N}$ be the corresponding null event.
	Then $\inf_{\|\vs-\vt\|<\epsilon} a_\vs = \inf_{\|\vs-\vt\|<\epsilon, \vs \in\mathcal{G}} a_\vs$ is measurable as a countable infimum, such that $a^*_\vt$ is a random variable.
	The lower semicontinuity of $a^*_\vt$ holds since
	\begin{align*}
		a^*_\vt = \lim_{\epsilon\to 0} \inf_{\|\vs-\vt\|<\epsilon} a_\vs \geq \lim_{\epsilon\to 0} \inf_{\|\vs-\vt\|<\epsilon} a^*_\vs \geq \lim_{\epsilon\to 0} \inf_{\|\vs-\vt\|<2\epsilon} a_\vs = a^*_\vt.
	\end{align*}
	
	To establish separability, fix some $\vt\in [0,T]^d$. 
	Since $a^*_t = \lim_{\epsilon\to 0} \inf_{\|\vs-\vt\|<\epsilon, \vs\in\mathcal{G}} a_\vs$, there exists a sequence $\vt_n\in \mathcal{G}$ such that $\vt_n\to \vt$ and $a_{\vt_n}\to a^*_\vt$ as $\to\infty$.
	Then, since $a^*_{\vt_n}\leq a_{\vt_n}$,
	\begin{align*}
		\limsup_{n\to\infty} a^*_{\vt_n} \leq  \limsup_{n\to\infty} a_{\vt_n} = a^*_{\vt}.
	\end{align*}
	By the lower semicontinuity of $a^*_\vt$, we also have $\liminf a^*_{\vt_n}\geq a^*_{\vt}$, such that $a^*_{\vt_n}\to a^*_{\vt}$ as $n\to\infty$.
	Hence, $a^*_\vt$ is separable with the same countable dense set $\mathcal{G}$ and the same null event $\mathcal{N}$ as $a_\vt$.	
\end{proof}

\begin{proof}[Proof of Theorem \ref{thm:KC}]
	The proof mimics central arguments of the standard Kolmogorov-Centsov theorem (see e.g.\ \cite[Thm. 2.2.8, Prob. 2.2.9]{Karatzas1998}). For completeness, we repeat some arguments given there.
	Fix some value $\epsilon>0$ such that \eqref{eqn:moment-condition} holds, and choose some $\gamma>0$. 
	Let $T=1$ without loss of generality. 
	Furthermore, conditions and claims of the Theorem are unaffected by choosing an equivalent norm, such that we may assume that $\|\cdot\|=\|\cdot\|_\infty$.
	Denote by $\mathcal{D}_n$ the discrete grid of mesh size $2^{-n}$,  
	\begin{align*}
		\mathcal{D}_n = \{ (k_1 2^{-n},\ldots, k_d 2^{-n}) : k_i=1,\ldots, 2^{n}-1, i=1,\ldots, d \}\subset(0,1)^d.
	\end{align*}
	For $n$ such that $2^{-n}<\epsilon$, Markov's inequality and the union bound yield
	\begin{align*}
		&\quad P\left( \bigcup_{\vs,\vt\in \mathcal{D}_n, \|\vs-\vt\|_\infty=2^{-n}} \left\{ |Y_{\vs} - Y_{\vt}| \geq 2^{-n\left(\underline{a}_{[\vs,\vt]}-\gamma\right)} \right\} \right) \\
		&\leq \sum_{\vs,\vt\in \mathcal{D}_n, \|\vs-\vt\|_\infty=2^{-n}} P\left( |Y_{\vs} - Y_{\vt}| \geq 2^{-n\left(\underline{a}_{[\vs,\vt]}-\gamma\right)}  \right) \\
		&\leq 3^d|\mathcal{D}_n| C(p,\epsilon) 2^{-nd} 2^{-n\gamma p} \\
		&= 3^dC(p,\epsilon) 2^{-n\gamma p}.
	\end{align*}
	Note furthermore that $\sum_{n=1}^\infty 2^{-n\gamma p}<\infty$. 
	Thus, the Borel-Cantelli Lemma yields that almost surely, there exists a (random) $N$ such that for all $n\geq N$ and all $\vs,\vt\in\mathcal{D}_n$, $\|\vs-\vt\|_\infty=2^{-n}$, 
	\begin{align}
		\left|Y_{\vs} - Y_{\vt}\right| < 2^{-n\left( \underline{a}_{[\vs,\vt]}-\gamma\right)}. \label{eqn:dyadic-bound-1}
	\end{align}
	We denote this event by $\tilde{\Omega}$, $P(\tilde{\Omega})=1$. 
	Denote by $\mathcal{D} = \bigcup_{n\in\N} \mathcal{D}_n$ the set of dyadic rationals in $(0,1)^d$. 	
	For two values $\vs,\vt\in\mathcal{D}_n$, let $l(\vs,\vt)$ be the unique integer $l$ such that $\|\vs-\vt\|_\infty \in [2^{-l-1}, 2^{-l})$.
	Using the construction of \cite[p.\ 341]{Potthoff2009}, there exist sequences ${\vs}^k, {\vt}^k \in \mathcal{D}_k, k=l(\vs,\vt),\ldots, n$, such that 
	\begin{align*}
		{\vs}^n=\vs,\quad \vt^n=\vt, \quad \vs^{l(\vs,\vt)} = \vt^{l(\vs,\vt)}, \\
		\|\vs^k - \vt^k\|_\infty \leq \|\vs^{k+1}-\vt^{k+1}\|_\infty, \\
		\|{\vs}^k - {\vs}^{k+1}\|_\infty= 2^{-(k+1)},\quad  \|\vt^k-\vt^{k+1}\|_\infty= 2^{-(k+1)} \;.
	\end{align*}
	In particular, $\|\vs^{k} - \vs\|_\infty \leq \sum_{j=l+1}^{n} 2^{-j} \leq 2^{-l} \leq 2\|\vs - \vt\|_\infty $ for $k=l(\vs,\vt),\ldots,n$, and $\|\vt^{k} - \vt\|_\infty \leq 2\|\vs - \vt\|_\infty$. 
	
	We denote by $\mathcal{B}$ the collection of closed sets in $[0,1]^d$ and introduce the notation $\underline{a}^\epsilon_{B}:= \inf_{\vr \in B^{\epsilon} \cap [0,1]^d} a_\vr$ for any closed set $B \in \mathcal{B}$ and $\epsilon>0$, where $B^{\epsilon} = \bigcup_{\vt\in B} B^{\|\cdot\|_\infty}_{\epsilon}(\vt)$ denotes the open $\epsilon$-neighborhood around $B$.
	Whenever $\|\vs-\vt\|_\infty < \epsilon$, we thus have that $\vs^k, \vt^k\in [\vs,\vt]^\epsilon$, such that $\underline{a}_{[\vs^k, \vt^k]} \geq \underline{a}^\epsilon_{[\vs,\vt]}$.
	If $2^{-N-1}\leq\|\vs-\vt\|_\infty < 2^{-N}$, then iterating \eqref{eqn:dyadic-bound-1} yields the inequality 
	\begin{align}
		\left|Y_{\vs}-Y_{\vt}\right| &\leq \sum_{k=N+1}^n |Y_{{\vs}^k}-Y_{{\vs}^{k+1}}| + |Y_{\vt^k}-Y_{\vt^{k+1}}| \nonumber \\
		&\leq \sum_{k=N+1}^\infty 2^{-k(\underline{a}_{[\vs,\vt]}^\epsilon - \gamma)} \nonumber \\
		&\leq 2^{-N (\underline{a}_{[\vs,\vt]}^\epsilon-\gamma)} \nonumber \\
		&\leq \left( 2 \|\vs-\vt\|_\infty\right)^{\underline{a}_{[\vs,\vt]}^\epsilon-\gamma }. \nonumber
	\end{align}
	 Introduce the event $E_\gamma = \{\gamma< \inf_{\vt\in (0,1)^d} a_\vt\}\cap \tilde{\Omega}$ such that $P(E_\gamma^c)\to 0$ as $\gamma\downarrow 0$.
	On the event $E_\gamma$, $Y_{\vt}$ is uniformly continuous on the set $\mathcal{D}$. 
	We define $\tilde{Y}_{\vt}=Y_{\vt}$ for $\vt\in \mathcal{D}$, and $\tilde{Y}_{\vt} = \lim_{\vs\to \vt, \vs\in\mathcal{D}} Y_{\vs}$ for $\vt\in [0,1]^d\setminus \mathcal{D}$. 
	The limit is independent of the specific sequence $\vs_n\to \vt$ due to the uniform continuity of $Y_{\vt}$.
	By monotonicity of the events $E_\gamma$, we may also define $\tilde{Y}_\mathbf{t}$ on the event $E_0 = \cup_{\gamma>0} E_\gamma$, with $P(E_0)=1$.
	Set $\tilde{Y}_\vs \equiv 0$ on $E_0^c$ without loss of generality, such that $\tilde{Y}_s$ is continuous for all $\omega\in\Omega$. 
	Now note that \eqref{eqn:moment-condition} yields in particular that $Y_{\vt}$ is continuous in probability, such that for each $\vt\notin\mathcal{D}$, $\tilde{Y}_{\vt} = Y_{\vt}$ almost surely.
	That is, $\tilde{Y}_{\vt}$ is a uniformly continuous modification of $Y_{\vt}$.
	For this modification, we have
	\begin{align*}
		P\left( \limsup_{m\to\infty} \sup_{\vs,\vt\in\mathcal{D}, \|\vs-\vt\|<\epsilon_m} \frac{|\tilde{Y}_{\vs} - \tilde{Y}_{\vt}|}{ \|\vs-\vt\|^{\underline{a}_{[\vs,\vt]}^{\epsilon_m}-\gamma}} \leq 2  \right) \geq P(E_\gamma),\quad \forall \gamma>0. 
	\end{align*}
	Here and in the following, $\epsilon_m = 2^{-m}$, where we use a countable sequence to ensure measurability.
	Then 
	\begin{align*}
		&\quad P\left( \forall \gamma>0\;: \limsup_{m\to\infty} \sup_{\vs,\vt\in \mathcal{D}, \|\vs-\vt\|<\epsilon_m} \frac{|\tilde{Y}_{\vs} - \tilde{Y}_{\vt}|}{ \|\vs-\vt\|^{\underline{a}^{\epsilon_m}_{[\vs,\vt]}-\gamma}} \leq 2 \right) \\
		&= \lim_{\gamma\to 0} P\left( \limsup_{m\to\infty} \sup_{\vs,\vt\in\mathcal{D}, \|\vs-\vt\|<\epsilon_m} \frac{|\tilde{Y}_{\vs} - \tilde{Y}_{\vt}|}{ \|\vs-\vt\|^{\underline{a}_{[\vs,\vt]}^{\epsilon_m}-\gamma}} \leq 2  \right) \\
		& \geq \lim_{\gamma\to 0} P(E_\gamma)= 1.
	\end{align*}
	
	By virtue of the continuity of $s\mapsto \tilde{Y}_s$, for all $\epsilon>0$ and $\delta\in\R$, and all $\omega\in\Omega$,
	\begin{align}
		\sup_{\vs,\vt\in [0,1]^d,\, \|s-t\| \geq \epsilon} \frac{|\tilde{Y}_s-\tilde{Y}_t|}{\|\vs-\vt\|^{\delta}} < \infty. \label{eqn:bounded-compact}
	\end{align}
	Hence,
	\begin{align*}
		\tilde{E}_0 &= \left\{\forall \gamma>0\;: \limsup_{m\to\infty} \sup_{\vs,\vt\in \mathcal{D}, \|\vs-\vt\|<\epsilon_m} \frac{|\tilde{Y}_{\vs} - \tilde{Y}_{\vt}|}{ \|\vs-\vt\|^{\underline{a}^{\epsilon_m}_{[\vs,\vt]}-\gamma}} <\infty \right\}  \\
		&= \left\{ \forall \gamma>0\;\forall M\in\N\;\exists m\geq M\;\forall B\in\mathcal{B}\;:  \sup_{\vs,\vt\in B\cap\mathcal{D}, \|\vs-\vt\|<\epsilon_m} \frac{|\tilde{Y}_{\vs} - \tilde{Y}_{\vt}|}{ \|\vs-\vt\|^{\underline{a}^{\epsilon_m}_{[\vs,\vt]}-\gamma}} <\infty \right\}  \\
		&\subset \left\{ \forall \gamma>0\;\forall M\in\N\;\exists m\geq M\;\forall B\in\mathcal{B}\;: \sup_{\vs,\vt\in B\cap\mathcal{D}, \|\vs-\vt\|<\epsilon_m} \frac{|\tilde{Y}_{\vs} - \tilde{Y}_{\vt}|}{ \|\vs-\vt\|^{\underline{a}^{2\epsilon_m}_{B}-\gamma}} <\infty \right\}   \\
		&\subset \left\{ \forall \gamma>0\;\forall M\in\N\;\exists m\geq M\;\forall B\in\mathcal{B}\;:  \sup_{\vs,\vt\in B\cap\mathcal{D}} \frac{|\tilde{Y}_{\vs} - \tilde{Y}_{\vt}|}{ \|\vs-\vt\|^{\underline{a}^{2\epsilon_m}_{B}-\gamma}} <\infty \right\}   \\
		\intertext{where we used \eqref{eqn:bounded-compact} and the fact that $\underline{a}_B^{2\epsilon} \leq \underline{a}_{[\vs,\vt]}^{\epsilon}$ for $\vs,\vt\in B$ with $\|\vs-\vt\|< \epsilon$,}
		&\subset \left\{ \forall \gamma>0\;\forall B\in\mathcal{B}\;\forall M\in\N\;\exists m\geq M\;:  \sup_{\vs,\vt\in B\cap\mathcal{D}} \frac{|\tilde{Y}_{\vs} - \tilde{Y}_{\vt}|}{ \|\vs-\vt\|^{\underline{a}^{\epsilon_m}_{B}-\gamma}} <\infty \right\}   \\
		&= \left\{ \forall \gamma>0\;\forall B\in\mathcal{B}\;\forall M\in\N\;\exists m\geq M\;:  \sup_{\vs,\vt\in B} \frac{|\tilde{Y}_{\vs} - \tilde{Y}_{\vt}|}{ \|\vs-\vt\|^{\underline{a}^{\epsilon_m}_{B}-\gamma}} <\infty \right\},
		\intertext{since $\tilde{Y}_\vs$ is uniformly continuous,}
		&\subset \left\{ \forall \gamma>0\;\forall B\in\mathcal{B}:  \sup_{\vs,\vt\in B} \frac{|\tilde{Y}_{\vs} - \tilde{Y}_{\vt}|}{ \|\vs-\vt\|^{\underline{a}_{B}-2\gamma}} <\infty \right\}  \\
		&= \left\{ \forall \gamma>0\;\forall B\in\mathcal{B}:  \sup_{\vs,\vt\in B} \frac{|\tilde{Y}_{\vs} - \tilde{Y}_{\vt}|}{ \|\vs-\vt\|^{\underline{a}_{B}-\gamma}} <\infty \right\}
	\end{align*}
	In the last step, we use that $\underline{a}^\epsilon_B\to \underline{a}_B$ as $\epsilon\to 0$ by virtue of Lemma \ref{lem:LSC}, such that $|\underline{a}_B - \underline{a}^\epsilon_B|< \gamma$ for $\epsilon$ sufficiently small.
	
	Hence, we have shown that
	\begin{align}
		\tilde{E}_0 \subset \left\{ \forall \gamma>0\;\forall B\in\mathcal{B}:  \sup_{\vs,\vt\in B} \frac{|\tilde{Y}_{\vs} - \tilde{Y}_{\vt}|}{ \|\vs-\vt\|^{\underline{a}_{B}-\gamma}} <\infty \right\}. \label{eqn:KC-subset}
	\end{align}
	The set $\tilde{E}_0$ is measurable and satisfies $P(\tilde{E}_0)=1$.
	However, for a general stochastic process $\tilde{Y}_\vs$, the right hand side in \eqref{eqn:KC-subset} is not necessarily measurable. 
	To fix this, we set $\tilde{Y}_\vs = 0$ on $\tilde{E}_0^c$, such that $\tilde{Y}_\vs$ is still a modification of $Y_\vs$.
	This modification satisfies $\left\{ \forall \gamma>0\;\forall B\in\mathcal{B}:  \sup_{\vs,\vt\in B} \frac{|\tilde{Y}_{\vs} - \tilde{Y}_{\vt}|}{ \|\vs-\vt\|^{\underline{a}_{B}-\gamma}} <\infty \right\} = \Omega$, such that the latter event is measurable and has probability $1$.
	
\end{proof}

\begin{proof}[Proof of Corollary \ref{cor:KC}]
	For any $\delta>0$, choose $p$ large enough such that $\delta p /2 \geq d$.
	Introduce the event $E_\delta = \{\inf_{\vt\in (0,T)^d} a_\vt >\delta\}\in \mathcal{F}$, such that $P(E_\delta)\uparrow 1$ as $\delta\downarrow 0$.	
	Then, for $\|\vs-\vt\|\leq \epsilon$ and $\epsilon$ sufficiently small, 
	\begin{align}
		\E  \left| \frac{Y_{\vs}-Y_\vt}{\|\vs-\vt\|^{\underline{a}_{[\vs,\vt]}-\delta}} \right|^p \mathds{1}_{E_\delta} 
		\;=\; \E  \left| \frac{Y_{\vs}-Y_\vt}{\|\vs-\vt\|^{\underline{a}_{[\vs,\vt]}-\frac{\delta}{2}}} \right|^p \mathds{1}_{E_\delta}  \|\vs-\vt\|^{\frac{\delta}{2} p} 
		\;\leq\; C(p, \epsilon, \delta/2)  \|\vs-\vt\|^d \label{eqn:cor-1}
	\end{align}
	We may now apply Theorem \ref{thm:KC} for the random fields $Y_\vt \mathds{1}_{E_\delta} $ and $\tilde{a}_\vt = (a_\vt-\frac{\delta}{2})\mathds{1}_{E_\delta} + \frac{1}{2} \mathds{1}_{E_\delta^c}$.
	That is, there exists a modification $\tilde{Y}^\delta_\vt$ of $Y_\vt\mathds{1}_{E_\delta} $ such that, for all $\delta\in(0,a)$, 
	\begin{align*}
		1 &= P\left( \forall \gamma>0\;\forall B\in\mathcal{B}:  \sup_{\vs,\vt\in B, \vs\neq\vt} \frac{|\tilde{Y}^\delta_{\vs} - \tilde{Y}^\delta_{\vt}|}{ \|\vs-\vt\|^{\underline{a}_B-\frac{\delta}{2}-\gamma}} <\infty \right) \\
		&= P\left( \forall \gamma>0\;\forall B\in\mathcal{B}:  \sup_{\vs,\vt\in B\cap \Q^d, \vs\neq\vt} \frac{|\tilde{Y}^\delta_{\vs} - \tilde{Y}^\delta_{\vt}|}{ \|\vs-\vt\|^{\underline{a}_B-\frac{\delta}{2}-\gamma}} <\infty \right).
	\end{align*}
	With probability one, it holds that $\tilde{Y}_t^\delta = Y_t \mathds{1}_{E_\delta}$ for all $\vt \in [0,T]^d \cap \Q^d$, such that
	\begin{align*}
		1	&= P\left( \forall \gamma>0\;\forall B\in\mathcal{B}:  \sup_{\vs,\vt\in B\cap \Q^d, \vs\neq\vt} \frac{|Y_{\vs} - Y_{\vt}|}{ \|\vs-\vt\|^{\underline{a}_B-\frac{\delta}{2}-\gamma}} \mathds{1}_{E_\delta} <\infty \right).
	\end{align*}
	In particular,
	\begin{align*}
		P\left( \forall \gamma>0\;\forall B\in\mathcal{B}:  \sup_{\vs,\vt\in B\cap \Q^d, \vs\neq\vt} \frac{|Y_{\vs} - Y_{\vt}|}{ \|\vs-\vt\|^{\underline{a}_B-\frac{\delta}{2}-\gamma}}  <\infty \right) \geq P(E_\delta).
	\end{align*}
	The event on the left hand side is decreasing in $\delta\downarrow 0$.
	Hence,
	\begin{align*}
		1 = \lim_{\delta\downarrow 0} P(E_\delta)&\leq \lim_{\delta\downarrow 0} P\left(\forall \gamma>0\;\forall B\in\mathcal{B}:  \sup_{\vs,\vt\in B\cap \Q^d} \frac{|Y_{\vs} - Y_{\vt}|}{ \|\vs-\vt\|^{\underline{a}_B-\frac{\delta}{2}-\gamma}} <\infty \right) \\
		&= P\left( \forall \delta>0\; \forall \gamma>0\;\forall B\in\mathcal{B}:  \sup_{\vs,\vt\in B\cap \Q^d} \frac{|Y_{\vs} - Y_{\vt}|}{ \|\vs-\vt\|^{\underline{a}_B-\frac{\delta}{2}-\gamma}} <\infty \right) \\
		&= P\left( \forall \gamma>0\;\forall B\in\mathcal{B}:  \sup_{\vs,\vt\in B\cap \Q^d} \frac{|Y_{\vs} - Y_{\vt}|}{ \|\vs-\vt\|^{\underline{a}_B-\gamma}} <\infty \right).
	\end{align*}
	In particular, $Y_\vt$ is almost surely uniformly continuous on $B \cap \Q^d$ for any closed set $B\subset[0,T]^d$. 
	Hence, there exists a continuous process $\tilde{Y}_\vt$ defined on $[0,T]^d$ such that $\tilde{Y}_\vt=Y_\vt$ for all $\vt\in[0,T]^d\cap \Q^d$, and
	\begin{align*}
		P\left( \forall \gamma>0\;\forall B\in\mathcal{B}:  \sup_{\vs,\vt\in B} \frac{|\tilde{Y}_{\vs} - \tilde{Y}_{\vt}|}{ \|\vs-\vt\|^{\underline{a}_B-\gamma}} <\infty \right) = 1.	
	\end{align*}
	By \eqref{eqn:cor-1}, $Y_\vt\mathds{1}_{E_\delta}$ is continuous in probability.
	Hence, $Y_\vt$ is also continuous in probability since
	\begin{align*}
		P\left( |Y_\vs-Y_\vt| > \epsilon \right) \leq P\left( |Y_\vs-Y_\vt|\mathds{1}_{E_\delta} > \epsilon \right) + P(E_\delta^c),
	\end{align*}
	and since $P(E_\delta^c)\to 0$ as $\delta\to 0$.
	In particular, we find that $Y_\vt=\tilde{Y}_\vt$ almost surely for all $\vt\in[0,T]^d$, i.e.\ $\tilde{Y}_\vt$ is a modification of $Y_\vt$.
	This establishes \eqref{eqn:KC-result}. 
\end{proof}

\subsection*{Proofs of Section \ref{sec:mbm}}

To apply the results of Section \ref{sec:KC}, we need to bound higher moments of stochastic processes defined via stochastic integrals.
The main technical tool applied in the derivations below is the Burkholder-Davis-Gundy (BDG) inequality, see e.g.\ \cite[Theorem IV.48]{Protter2005}. 
For a standard Wiener process $W_t$ and an adapted, left-continuous integrand $h_s$, and any $p\geq 1$, it yields that
\begin{align}
	\E \left| \int_{-\infty}^{\infty} h_s \, dW_s\right|^p \leq C_p \E \left| \int_{-\infty}^{\infty} h_s^2\, ds \right|^\frac{p}{2}, \label{eqn:BDG}
\end{align}
for a universal constant $C_p$.
In the following, it will thus typically suffice to derive bounds on various integral expressions corresponding to the right hand side of \eqref{eqn:BDG}.

\begin{proof}[Proof of Lemma \ref{lem:inc-moment}]
	Let $\epsilon \in (0,1/2)$, $\delta=\omega(2 \epsilon)$, and $h\in(0,\epsilon)$. 
	We decompose
	\begin{align*}
		X_{t+h} - X_t &= \int_{-\infty}^{t-\epsilon} \left[g_s(t+h)-g_s(t) \right] dW_s +  \int_{t-\epsilon}^{t} \left[g_s(t+h)-g_s(t) \right] dW_s \\
		&\qquad + \int_{t}^{t+h} g_s(t+h) dW_s \\
		&= D_\epsilon + E_\epsilon + F.
	\end{align*}
	Note that $\E |D_\epsilon + E_\epsilon + F|^p \leq 3^p(\E |D_\epsilon|^p + \E|E_\epsilon|^p + \E|F|^p)$, such that it suffices to bound the three terms individually. 
	In the following we denote $\underline{H}_t = \inf_{s\in[t-\epsilon,t+\epsilon]} H_s$. 
	We consider the continuous case first, i.e.\ \eqref{eqn:inc-moment-res-1}.
	
	\underline{Bounding $D_\epsilon$:}
	To bound the term $D_\epsilon$, we apply the mean value theorem to obtain for some $\tilde{h}_{s,t} \in[0,h]$, 
	\begin{align*}
		&\quad \int_{-\infty}^{t-\epsilon} |g_s(t+h) - g_s(t)|^2\, ds \\
		&= h^2 \int_{-\infty}^{t-\epsilon} \left| \partial_t g_s(t+\tilde{h}_{s,t}) \right|^2\,ds \\
		&\leq \overline{L}^2 h^2 \left[\int_{-\infty}^{t-\epsilon} |t+\tilde{h}_{s,t}-s|^{2H_s-3} \mathds{1}_{|t+\tilde{h}_{s,t}-s|\leq 1} \;+\; |t+\tilde{h}_{s,t} - s|^{-2R_s} \mathds{1}_{|t+\tilde{h}_{s,t} - s|>1}\,ds \right]\\
		&\leq \overline{L}^2 h^2 \left[\int_{-\infty}^{t-\epsilon} |t-s|^{2\underline{H}-3} \mathds{1}_{|t+\tilde{h}_{s,t}-s|\leq 1} \;+\; |t - s|^{-2\underline{R}} \mathds{1}_{|t+\tilde{h}_{s,t} - s|>1}\,ds \right]\\
		&\leq \overline{L}^2 h^2 \left[\int_{t-1}^{t-\epsilon} |t-s|^{2\underline{H}-3} \, ds \;+\; \int_{-\infty}^{t-\frac{1}{2}}|t - s|^{-2\underline{R}} \,ds \right]\\		
		&= \overline{L}^2 h^2 \left[\int_{\epsilon}^{1} z^{2\underline{H}-3}\, dz +  \int_{\frac{1}{2}}^{\infty} |z|^{-2\underline{R}}\, dz \right]\\
		&\leq 2\overline{L}^2h^2 \left[ \frac{\epsilon^{2\underline{H}-2}}{2-2\underline{H}} + \frac{2^{2\underline{R}-1}}{2\underline{R}-1} \right]\\
		&\leq C(\underline{H},\overline{L}, \underline{R}) h^2 \epsilon^{2\underline{H}-2}.
	\end{align*}
	Then, the Burkholder-Davis-Gundy inequality yields for any $p>1$,
	\begin{align}
		\E \left|D_\epsilon h^{-H_t} \right|^p 
		&\leq h^{-p\overline{H}} \E \left|D_\epsilon  \right|^p \label{eqn:Deps-bound}\\
		&\leq h^{-p\overline{H}} C_p \E \left[ \int_{-\infty}^{t-\epsilon} |g_s(t+h) - g_s(t)|^2\, ds \right]^\frac{p}{2} \nonumber\\
		&\leq h^{p(1-\overline{H})} C(\overline{H},\overline{L},\underline{R},p) \epsilon^{p(\underline{H}-1)} \nonumber\\
		&\leq C(\overline{H},\overline{L},\underline{R},p, \epsilon), \nonumber
	\end{align}
	where $C(\overline{H},\overline{L},\underline{R},p)$ is a constant depending on $\overline{H}$,$\overline{L}$, $\underline{R}$, and $p$.
	
	\underline{Bounding $E_\epsilon$:}
	Special care is needed if $H_s$ is close to $\frac{1}{2}$.
	For this reason, we define for any $\delta>0$ the process \begin{align*}
		H_s^\delta = \begin{cases}
			\frac{1}{2}-\delta,& H_s \in \left[\frac{1}{2}-\delta, \frac{1}{2}+\delta \right],\\
			H_s, & \text{otherwise}.
		\end{cases}
	\end{align*}
	We also extend the notation $\underline{H}_t^\delta = \inf_{r\in[t-\epsilon,t+\epsilon]} H_r^\delta$.
	Now write $g_s(t+h) - g_s(t)= \int_t^{t+h} g'_s(r)\, dr$.
	Since $H^\delta_s \leq H_s$, we obtain
	\begin{align}
		\int_{t-\epsilon}^{t} |g_s(t+h) - g_s(t)|^2\, ds 
		&\leq \int_{t-\epsilon}^{t} \left| \int_t^{t+h} |g'_s(r)|\, dr \right|^2\, ds \nonumber\\
		&\leq \int_{t-\epsilon}^{t} \left| \int_t^{t+h} L_s |r-s|^{H_s - \frac32}  \, dr \right|^2\, ds \nonumber\\
		&\leq \int_{t-\epsilon}^{t} \left| \int_t^{t+h} L_s |r-s|^{\underline{H}_t^\delta - \frac32}  \, dr \right|^2\, ds \nonumber\\
		&\leq \frac{\overline{L}^2}{\left|\underline{H}_t^\delta - \frac{1}{2}\right|^2}\int_{t-\epsilon}^t  \left| |t+h-s|^{\underline{H}_t^\delta-\frac{1}{2}} - |t-s|^{\underline{H}_t^\delta-\frac{1}{2}} \right|^2\, ds \nonumber\\
		&\leq \overline{L}^2  \frac{h^{2\underline{H}_t^\delta}}{\delta^2} \int_{0}^{\infty} \left|  (r+1)^{\underline{H}_t^\delta-\frac{1}{2}} - r^{\underline{H}_t^\delta - \frac{1}{2}}\right|^2\, dr. \label{eqn:Eeps-bound}
	\end{align}
	The latter equality follows by substitution. 
	Note that $|\underline{H}_t^\delta-1/2| \geq \delta$.
	Moreover, there exists a constant $C(\underline{H}, \overline{H})$ such that for all $H\in(\underline{H},\overline{H})$,
	\begin{align}
		&\quad \int_{0}^{\infty}  \left|  (r+1)^{H-\frac{1}{2}} - r^{H - \frac{1}{2}}\right|^2\, dr  \nonumber\\
		&= \int_{0}^{1} \left|  (r+1)^{H-\frac{1}{2}} - r^{H - \frac{1}{2}}\right|^2\, dr + \int_{1}^{\infty}  \left|  (r+1)^{H-\frac{1}{2}} - r^{H - \frac{1}{2}}\right|^2\, dr \nonumber\\
		&\leq 2\int_{0}^{1} \left[(r+1)^{2 \overline{H}-1} + r^{2 \underline{H} - 1}\right]\, dr + 2\int_1^\infty \sup_{H\in(\underline{H}, \overline{H})}|2r|^{2H-3} \,dr \nonumber\\
		& \leq 2 \int_1^2 r^{2\overline{H}-1}\, dr + 2\int_0^1 r^{2\underline{H}-1}\, dr + 2 \int_1^\infty r^{2\overline{H}-3}\, dr 
		\qquad \leq\qquad C(\underline{H}, \overline{H}). \label{eqn:subst-bound}
	\end{align}	
	We now apply the Burkolder-Davis-Gundy inequality to the local martingale defined by $h^{-H_{t-\epsilon}} \int_{t-\epsilon}^r \left[g_s(t+h) - g_s(t)\right]\, dW_s$.
	By additionally exploiting the continuity of $H_t$, we find that
	\begin{align*}
		\E | E_\epsilon h^{-H_t}|^p 
		& \leq h^{-p\,\omega(\epsilon)} \E | h^{-H_{t-\epsilon}} E_\epsilon|^p \\
		&\leq  h^{-p\,\omega(\epsilon)} C_p \E \left[ h^{-2H_{t-\epsilon}} \int_{t-\epsilon}^t |g_s(t+h)-g_s(t)|^2\,ds \right]^\frac{p}{2} \\
		&\leq h^{-p\, \omega(2\epsilon)} C_p \E \left[ \overline{L}^2 \delta^{-2} C(\underline{H},\overline{H}) h^{2\underline{H}_t^\delta} h^{-2H_{t-\epsilon}} \right]^\frac{p}{2} \\
		&\leq h^{-2p\, \omega(2\epsilon)} C_p \E \left[ \overline{L}^2 \delta^{-2} C(\underline{H},\overline{H}) h^{2(\underline{H}_t^\delta - \underline{H}_t)} \right]^\frac{p}{2} \\
		&\leq h^{-2p\, \omega(2\epsilon)-2p\delta} C(\overline{L}, \underline{H},\overline{H}, p, \delta),
	\end{align*}
	since $|\underline{H}_t^\delta - \underline{H}_t| \leq \delta$.
	Here, it is crucial to use the value $H_{t-\epsilon}$ instead of $H_t$ in the expectation to ensure that the integrand is adapted, such that the BDG inequality is applicable.
	
	\underline{Bounding $F$:} 
	Since $g_s(t)=0$ for $s>t$, we obtain from \eqref{eqn:cond-B-3} 
	\begin{align}
		\int_{t}^{t+h} |g_s(t+h)|^2\, ds  
		\quad\leq\quad \overline{L} \int_{0}^{h} r^{2 \underline{H}_t -1}\, dr 
		\quad\leq\quad C(\underline{H}, \overline{L}) h^{2 \underline{H}_t} . \label{eqn:F-bound}
	\end{align}
	Then, just as for the term $E_\epsilon$, the Burkholder-Davis-Gundy inequality and the continuity of $H_t$ yields \begin{align*}
		\E |h^{-H_t} F|^p \leq h^{-2p\, \omega(2\epsilon)} C(\overline{L}, \underline{H},\overline{H}, p).
	\end{align*}
	This completes the proof of \eqref{eqn:inc-moment-res-1}.
	
	\underline{Discontinuous $H_t$:}
	If $H_t$ is discontinuous, \eqref{eqn:Deps-bound} still applies and yields $\E|D_\epsilon h^{-H_t}|^p \leq C(\overline{H},\overline{L}, \underline{R},\epsilon, p)$.
	Moreover, if $H_t$ is deterministic, we obtain from \eqref{eqn:F-bound} and \eqref{eqn:Eeps-bound},
	\begin{align*}
		\E \left| h^{-\underline{H}^\epsilon_t} (E_\epsilon+F) \right|^p 
		&\leq C(\underline{H},\overline{H},\overline{L},p) h^{-2\underline{H}^\epsilon_t } \E \left| \int_{t-\epsilon}^{t+h} |g_s(t+h-s) - g_s(t-s)|^2\, ds \right|^\frac{p}{2} \\
		&\leq C(\underline{H},\overline{H},\overline{L},p, \delta) h^{- p \delta}.
	\end{align*}
	This establishes \eqref{eqn:inc-moment-res-3}.
	
	Moreover, for discontinuous but possibly random $H_t$, \eqref{eqn:F-bound} and \eqref{eqn:Eeps-bound} yield $\E|h^{-\underline{H}} (E_\epsilon+F)|^p \leq C(\underline{H},\overline{H},\overline{L}, p, \delta) h^{-p\delta}$ and $\E|h^{-\underline{H}} D_\epsilon|^p \leq C(\overline{H},\overline{L}, \underline{R}, p, \epsilon)$, such that \eqref{eqn:inc-moment-res-2} holds.
\end{proof}

\begin{proof}[Proof of Theorem \ref{thm:regularity-X}]
	Property \eqref{eqn:inc-moment-res-2} implies that, if $|s-t|< \epsilon$,
	\begin{align*}
		\E \left| \frac{X_s-X_t}{|s-t|^{\underline{H} -\delta}} \right|^p \leq C(\underline{H},\overline{H}, \overline{L}, \underline{R}, \epsilon, p, \delta) .
	\end{align*}
	Hence, Corollary \ref{cor:KC} is applicable and yields the existence of a continuous modification with Hölder exponent $\underline{H}$, such that \eqref{eqn:smoothness-lower} holds.
	
	Moreover, \eqref{eqn:inc-moment-res-1} implies that, for $s<t, |s-t|<\epsilon$,
	\begin{align*}
		\E \left| \frac{X_s-X_t}{|s-t|^{\underline{H}_{[s,t]} -\delta - 2\omega(2\epsilon)}} \right|^p \leq C(\underline{H},\overline{H}, \overline{L}, \underline{R}, \epsilon, p, \delta)
	\end{align*}
	Recall the notation $\underline{H}_{[s,t]}=\inf_{r\in[s,t]} H_r$.
	Since $\omega(\epsilon)\to 0$ as $\epsilon\to 0$, we may apply Corollary \ref{cor:KC} to obtain, for a modification $X_t$, 
	\begin{align}
		1 &= P\left( \forall \gamma>0\;  \forall \,0<a<b<T: \sup_{s,t\in[a,b]} \frac{|X_s-X_t|}{|s-t|^{\underline{H}_{[a,b]}-\gamma}} < \infty \right). \label{eqn:Ito-mBm-uniform-Holder}
	\end{align}
	Now use that the the continuity of $t\mapsto H_t$ implies that $\underline{H}_{[t-|h|, t+|h|]}\to H_t$ for each $t$, as $h\to 0$.
	Hence,
	\begin{align*}
		1 &= P\left( \forall \gamma>0,\;  t\in(0,T): \limsup_{|h|\to 0} \sup_{s,r\in[t-|h|,t+|h|]} \frac{|X_s-X_r|}{|s-r|^{\underline{H}_{[t-|h|,t+|h|]}-\gamma}} <\infty \right) \\
		&= P\left( \forall \gamma>0,\;  t\in(0,T): \limsup_{|h|\to 0}  \frac{|X_{t+h}-X_t|}{|h|^{\underline{H}_{[t-|h|,t+|h|]}-\gamma}} <\infty \right) \\
		&= P\left( \forall \gamma>0,\;  t\in(0,T): \limsup_{|h|\to 0}  \frac{|X_{t+h}-X_t|}{|h|^{H_t-\gamma}} <\infty \right) \\
		&= P\left( \forall \gamma>0,\;  t\in(0,T): \limsup_{|h|\to 0}  \frac{|X_{t+h}-X_t|}{|h|^{H_t-\gamma}} =0 \right) \\
		&= P\left( \alpha_t(X) \geq H_t \quad\forall t\in(0,T)\right).
	\end{align*}
	This establishes \eqref{eqn:smoothness-local}.
	
	For the discontinuous case \eqref{eqn:smoothness-deterministic} we replace $\underline{H}^\epsilon_t = \inf_{r\in[t-\epsilon, t+\epsilon]} H_r$ by its lower semicontinuous version $\underline{H}_t^{\epsilon, *}\leq \underline{H}^\epsilon_t$, given by $\underline{H}_t^{\epsilon,*}= \lim_{\delta\to 0} \inf_{s\in[t-\delta, t+\delta]} \underline{H}^\epsilon_s$.
	Just as in Theorem \ref{thm:KC}, we also denote, for $0<s<t<T$,
	\begin{align*}
		\underline{H}_{[s,t]}^{\epsilon,*} = \inf_{r\in[s, t]} \underline{H}_r^{\epsilon,*}.
	\end{align*}
	Then \eqref{eqn:inc-moment-res-3} yields, for all $s,t\in(0,T)$,
	\begin{align*}
		\E \left| \frac{X_t-X_s}{|t-s|^{\underline{H}_{[s,t]}^{\epsilon,*} - \delta}} \right| \leq C(\underline{H}, \overline{H}, \overline{L}, \underline{R}, \epsilon, p, \delta)
	\end{align*}
	Hence, Theorem \ref{thm:KC} yields, for any $\epsilon\in(0,1)$,
	\begin{align*}
		1 
		&= P\left( \forall \gamma>0 \; \forall\, 0<a<b<T : \,  \sup_{s,t\in [a,b]} \frac{|X_s - X_t|}{ |s-t|^{\underline{H}_{[a,b]}^{\epsilon,*}-\gamma}} <\infty \right).
	\end{align*}
	We may now proceed just as for the proof of \eqref{eqn:smoothness-local}, since $\underline{H}_{[t-|h|, t+|h|]}^{\epsilon,*} \to \underline{H}_t^{\epsilon,*}$, to obtain
	\begin{align*}
		1 = P\left( \alpha_t(X) \geq  \underline{H}_t^{\epsilon,*}, \;\forall t\in(0,T)\right).
	\end{align*}
	Now note that $\underline{H}_t^{\epsilon,*} \uparrow H_t^*$ as $\epsilon \to 0$. 
	Since $\epsilon>0$ may be chosen arbitrarily small, we conclude that 
	\begin{align*}
		P\left( \alpha_t(X) \geq  H_t^*,\;\forall t\in(0,T)\right) = 1,
	\end{align*}
	establishing \eqref{eqn:smoothness-deterministic}.
\end{proof}

\begin{proof}[Proof of Theorem \ref{thm:rescaling}]
	We will first establish that, for $q\in(0,1)$ sufficiently small,
	\begin{align}
		\begin{split}
			&h^{-H_t}\left(X_{t+hr} - X_t\right)\\
			&= h^{-H_{t-h^q}}\int_{t-h^q}^{\infty} \sigma_{t-h^q}\left[(t+hr-s)_+^{H_{t-h^q}-\frac{1}{2}} - (t-s)_+^{H_{t-h^q}-\frac{1}{2}} \right] \, dW_s + o_P(1). 
		\end{split} \label{eqn:rescaling-0}
	\end{align}
	Here, $o_P(1)$ denotes a family $\varepsilon_h$ of random variables such that $\varepsilon_h\to 0$ in probability as $h\to 0$. 
	The asymptotic decomposition \eqref{eqn:rescaling-0} holds for any fixed $t>0$ and $r\in(a,b)$ and allows us to establish convergence of the finite dimensional distributions.
	We do not need the $o_P(1)$ term to vanish uniformly in $r\in (a,b)$, as tightness will be established separately in step (vi) below. 
	
	\underline{Step (i):} 
	Note that we established in \eqref{eqn:Deps-bound} that
	\begin{align}
		\E\left|h^{-H_t}\int_{-\infty}^{t-h^q} \left[ g_s(t+hr) - g_s(t) \right] \,dW_s \right|^2 \leq C( \overline{H}, \overline{L}, \underline{R}, 2) h^{2(1-\overline{H})} h^{-2\,q\,(1-\underline{H})}, \label{eqn:rescaling-2}
	\end{align}
	which tends to zero for $q$ sufficiently small.
	In particular, 
	\begin{align*}
		h^{-H_t}\left(X_{t+hr} - X_t\right) 
		&= h^{-H_t}\int_{t-h^q}^{\infty} \left[g_s(t+hr)- g_s(t)\right] \, dW_s + o_P(1) \\
		&= h^{-H_{t-h^q}}\int_{t-h^q}^{\infty} \left[g_s(t+hr)- g_s(t)\right] \, dW_s \; (1+o_P(1)) + o_P(1).
	\end{align*}
	The latter step is valid because by the mean value theorem 
	\[|h^{H_t-H_{t-h^q}}-1| \leq \exp(-\log(h) \omega(h^q)) \left|\omega(h^q)\log(h) \right|,\]
	and the assumptions on $\omega$ imply that 
	\begin{align*}
		\lim_{h\to 0} \omega(h^q)\log(h) = \lim_{h\to 0} \omega(h) \log(h^\frac{1}{q}) = 0.
	\end{align*} 
	
	\underline{Step (ii):}
	Let $r>0$ without loss of generality.
	If $r<0$, we may just exchange the roles of $t$ and $t+hr$ in the following derivations.
	Define the integral kernel 
	\begin{align*}
		\tilde{g}_s(t) = \left[g_s(t) - \sigma_s(t-s)_+^{H_s-\frac{1}{2}} \right] \mathds{1}_{s\geq t-h^q}.
	\end{align*}
	Now observe that condition \nameref{cond:B2} implies that $\tilde{g}_s(t)$ satisfies condition \nameref{cond:g_s(t)} with exponent $\tilde{H}_s = H_s+\rho$.
	Then the same arguments as in the proof of Lemma \ref{lem:inc-moment}, i.e.\ \eqref{eqn:Eeps-bound} and \eqref{eqn:F-bound}, yield, for any $\delta>0$,
	\begin{align*}
		\int_{t-h^q}^{t+hr} \left|\tilde{g}_s(t+hr) - \tilde{g}_s(t)\right|^2\, ds 
		&\leq h^{2(\inf_{s\in[t-h^q, t+h]}H_s +\rho-2\delta )} C(\overline{L}, \underline{H},\overline{H}, \delta, \rho) \\
		&\leq h^{2(H_{t-h^q} - \omega(h^q + h|r|)+\rho-2\delta )} C(\overline{L}, \underline{H},\overline{H}, \delta, \rho).
	\end{align*}
	Now choose $\delta<\rho/4$.
	Since $h^{\omega(h^q+h|r|)}=\exp( \log(h) \omega(h^q+h|r|))\to 1$, the latter bound is of order $h^{2H_{t-h^q}+\rho}$.
	Then the Burkholder-Davis-Gundy inequality yields 
	\begin{align*}
		h^{-H_{t-h^q}} \int_{t-h^q}^\infty \left[\tilde{g}_s(t+hr)-\tilde{g}_s(t)\right]\, dW_s \quad\pconv \quad 0\quad \text{ as }h\to 0,
	\end{align*}
	such that
	\begin{align*}
		&\qquad h^{-H_{t-h^q}}\left(X_{t+hr} - X_t\right) \\
		&= h^{-H_{t-h^q}}\int_{t-h^q}^{\infty} \sigma_s\left[(t+hr-s)_+^{H_s-\frac{1}{2}}- (t-s)_+^{H_s-\frac{1}{2}}\right] \sigma_s\, dW_s \; (1+o_P(1)) + o_P(1).
	\end{align*}
	
	\underline{Step (iii):}
	In the integrand, we approximate the exponent $H_s$ by $H_{t-h^q}$.
	Indeed, the technical Lemma \ref{lem:technical-1} below establishes that 
	\begin{align*}
		&\quad \int_{t-h^q}^{\infty} \left|\left[(t+hr-s)_+^{H_s-\frac{1}{2}} - (t-s)_+^{H_s-\frac{1}{2}}\right] - \left[(t+hr-s)_+^{H_{t-h^q}-\frac{1}{2}} - (t-s)_+^{H_{t-h^q}-\frac{1}{2}}\right]\right|^2 \, ds \\
		&= \int_{-h^q}^{hr} \left|\left[(hr-s)_+^{H_{t+s}-\frac{1}{2}} - (-s)_+^{H_{t+s}-\frac{1}{2}}\right] - \left[(hr-s)_+^{H_{t-h^q}-\frac{1}{2}} - (-s)_+^{H_{t-h^q}-\frac{1}{2}}\right]\right|^2 \, ds \\
		&\leq \sup_{s\in[t-h^q, t+hr]} |H_s-H_{t-h^q}|\, |hr|^{2H_{t-h^q}-2\omega(h^q+h|r|)} \log(h|r|)^2 C(\underline{H}, \overline{H}) \\
		&\leq \omega(h^q+h|r|)^2 |hr|^{2H_{t-h^q}-2\omega(h^q+h|r|)} \log(h|r|)^2 C(\underline{H}, \overline{H}) \\
		&\leq \omega(2h^q)^2 h^{2H_{t-h^q}} \log(h)^2 C(\underline{H},\overline{H},r) \stackrel{P}{\to} 0
	\end{align*}
	because $\omega(2h^q)\log(h) \to 0$. 
	Since $\sigma_s$ is bounded, the Burkholder-Davis-Gundy inequality yields
	\begin{align*}
		&\qquad h^{-H_t}(X_{t+hr}-X_t) \\
		&= h^{-H_{t-h^q}} \int_{t-h^q}^{\infty}  \sigma_s\left[ (t+hr-s)_+^{H_{t-h^q}-\frac{1}{2}} - (t-s)_+^{H_{t-h^q}-\frac{1}{2}}\right]\, dW_s\; (1+o_P(1)) + o_P(1).
	\end{align*}
	
	\underline{Step (iv):}
	The term $\sigma_s$ in the integrand can be approximated by $\sigma_{t-h^q}$. 
	In fact, a substitution yields
	\begin{align*}
		&\quad h^{-2H_{t-h^q}}\int_{t-h^q}^{t\vee(t+hr)} |\sigma_s-\sigma_{t-h^q}|^2\left| \left[ (t+hr-s)_+^{H_{t-h^q}-\frac{1}{2}} - (t-s)_+^{H_{t-h^q}-\frac{1}{2}}\right] \right|^2\, ds \\
		&\leq \sup_{s\in[t-h^q, t]} |\sigma_s-\sigma_{t-h^q}|^2 \int_{-\infty}^{\infty} \left| (t+r-s)_+^{H_{t-h^q}-\frac{1}{2}} - (t-s)_+^{H_{t-h^q}-\frac{1}{2}}\right|^2\, ds.
	\end{align*}
	The second term is bounded by a constant $C(\overline{H}, \underline{H}, r)$, as established in \eqref{eqn:subst-bound}.
	Since $\sigma_s$ is continuous and bounded by a constant, we also have $\E (\sup_{s\in[t-h^q,t]} |\sigma_s-\sigma_{t-h^q}|^2)\to 0$ as $h\to 0$.
	Thus, 
	\begin{align*}
		&\qquad h^{-H_t} (X_{t+hr}-X_t) \\
		&= h^{-H_{t-h^q}} \int_{t-h^q}^\infty \sigma_{t-h^q}\left[ (t+hr-s)_+^{H_{t-h^q}-\frac{1}{2}} - (t-s)_+^{H_{t-h^q}-\frac{1}{2}}\right]\, dW_s \; (1+o_P(1)) + o_P(1),
	\end{align*}
	which establishes \eqref{eqn:rescaling-0}.
	
	\underline{Step (v):}
	To study the dependency with $Z$, let $Z_h = \E(Z|\F_{t-h^q})$. 
	Applying the martingale convergence theorem, it can be established that $Z_h \pconv \E(Z|\F_{t-})=Z$ as $h\to 0$.
	Hence, \eqref{eqn:rescaling-0} yields 
	\begin{align*}
		&\quad \begin{pmatrix}
			h^{-H_t} (X_{t+hr}-X_t) \\ Z 
		\end{pmatrix}
		\\ 
		&=
		\begin{pmatrix}
			h^{-H_{t-h^q}} \displaystyle\int_{t-h^q}^\infty \sigma_{t-h^q}\left[ (t+hr-s)_+^{H_{t-h^q}-\frac{1}{2}} - (t-s)_+^{H_{t-h^q}-\frac{1}{2}}  \right]\, dW_s \\
			Z_h
		\end{pmatrix} \; (1+o_P(1)) + o_P(1) \\
		&= \Xi_h \; (1+o_P(1)) + o_P(1).
	\end{align*}
	Note that both, $Z_h$ and the integrand of the stochastic integral, are independent of $W_s-W_{t-h^q}, s\geq t-h^q$. 
	Therefore, we may introduce a standard Brownian motion $\tilde{W}_t$ independent of the processes $W$ and $H$, such that
	\begin{align*}
		\Xi_h 
		&\deq \begin{pmatrix}
			h^{-H_{t-h^q}} \displaystyle\int_{t-h^q}^\infty \sigma_{t-h^q}\left[ (t+hr-s)_+^{H_{t-h^q}-\frac{1}{2}} - (t-s)_+^{H_{t-h^q}-\frac{1}{2}}  \right]\, d\tilde{W}_s \\
			Z_h
		\end{pmatrix}\\
		&\deq \begin{pmatrix}
			h^{-H_{t-h^q}} \displaystyle\int_{-\infty}^\infty \sigma_{t-h^q}\left[ (t+hr-s)_+^{H_{t-h^q}-\frac{1}{2}} - (t-s)_+^{H_{t-h^q}-\frac{1}{2}}  \right]\, d\tilde{W}_s \\
			Z
		\end{pmatrix} + o_P(1)\\
		&= \tilde{\Xi}_h+o_P(1).
	\end{align*}
	Here, we also applied \eqref{eqn:rescaling-2} in the second step.
	Now, self-similarity of $\tilde{W}_s$ and the dominated convergence theorem for stochastic integrals yield
	\begin{align*}
		\tilde{\Xi}_h 
		&\deq \begin{pmatrix}
			\displaystyle\int_{-\infty}^{\infty} \sigma_{t-h^q}\left[ (r-s)_+^{H_{t-h^q}-\frac{1}{2}} - (-s)_+^{H_{t-h^q} - \frac{1}{2}} \right]\, d\tilde{W}_s \\
			Z
		\end{pmatrix} \\ 
		&=\begin{pmatrix}
			\displaystyle\int_{-\infty}^{\infty} \sigma_t\left[ (r-s)_+^{H_t-\frac{1}{2}} - (-s)_+^{H_t - \frac{1}{2}} \right]\, d\tilde{W}_s \\
			Z
		\end{pmatrix} + o_P(1).
	\end{align*}
	
	The same reasoning applies for multiple values $r_1,\ldots, r_k\in\R$.
	Hence, we obtain the multivariate convergence in distribution by Slutsky, as $h\to 0$,
	\begin{align}
		\begin{pmatrix}
			h^{-H_t} \left(X_{t+hr_j}-X_t\right)_{j=1}^k
			\\ Z
		\end{pmatrix}
		\wconv 
		\begin{pmatrix}
			\left(\displaystyle\int_{-\infty}^{\infty} \sigma_t \left[(r_j-z)_+^{H_t-\frac{1}{2}} - (-z)_+^{H_t-\frac{1}{2}} \right]\, d\tilde{W}_z \right)_{j=1}^k \\ 
			Z
		\end{pmatrix}. \label{eqn:fidi-conv}
	\end{align}
	
	\underline{Step (vi):}
	Let $I=(a,b)$ be a finite interval with $a<0<b$. 
	To show that the weak convergence \eqref{eqn:fidi-conv} also holds in a functional sense, it suffices to show that the sequence of probability measures on $C(I)$ induced by $(\tilde{X}^{t,h}_r)_{r\in I}=h^{-H_t}(X_{t+hr}-X_t)_{r\in I}$ is tight.
	Then, for $r,s\in I$, $r<s$, $|r-s|<1$,
	\begin{align*}
		\E  \left|\frac{\tilde{X}^{t,h}_r-\tilde{X}^{t,h}_s}{|r-s|^{\underline{H}}}\right|^p 
		&= \E \left| \frac{X_{t+hr}-X_{t+hs}}{h^{H_t} |r-s|^{\underline{H}}} \right|^p \\
		&\leq \E \left| \frac{X_{t+hr}-X_{t+hs}}{h^{H_t-H_{t+hr}}|hr-hs|^{H_{t+hr}}} \right|^p \\
		&\leq \E \left| \frac{X_{t+hr}-X_{t+hs}}{|hr-hs|^{H_{t+hr}}} \right|^p h^{-p\omega(h|b-a|)}.
	\end{align*}
	Since $h^{\omega(h|b-a|)}\to 1$ by assumption, Lemma \ref{lem:inc-moment} yields for $h|r-s|\in(0,\epsilon)$, $\epsilon>0$ and $\delta>0$ sufficiently small,
	\begin{align*}
		\E  \left|\frac{\tilde{X}^{t,h}_r-\tilde{X}^{t,h}_s}{|r-s|^{\underline{H}}}\right|^p & \leq C(\underline{H}, \overline{H}, \overline{L}, \underline{R}, \epsilon, p,\delta)\, |hr-hs|^{-2p\omega(2\epsilon)-\delta}.
	\end{align*}
	Thus, Corollary \ref{cor:KC} applies and yields 
	\begin{align}
		P\left( \limsup_{\epsilon\to 0} \sup_{r,s\in I, |r-s|<\epsilon} \frac{|\tilde{X}^{t,h}_r - \tilde{X}^{t,h}_s|}{|r-s|^{\underline{H}/2}} \leq 1 \right) = 1. \label{eqn:rescaled-holder}
	\end{align}
	The fact that $\tilde{X}^{t,h}_0=0$ for all $h>0$, in combination with \eqref{eqn:rescaled-holder}, establishes tightness of $X^{t,h}_r$ in the space of continuous functions $C(I)$ \cite[Thm. 2.8.1]{Billingsley1968}.
\end{proof}

\begin{lemma}\label{lem:technical-1}
	Let $a_s, b_s\in [\underline{c},\overline{c}] \subset (0,1)$, and $|a_s-b_s|\leq \Delta$ for all $s\in\R$. 
	Then for all $h\in(0,1)$,
	\begin{align}
		\begin{split}
			&\qquad \int_{-\infty}^{0} \left| \left[(h-s)^{a_s-\frac{1}{2}} - (-s)^{a_s-\frac{1}{2}} \right] - \left[(h-s)^{b_s-\frac{1}{2}} - (-s)^{b_s-\frac{1}{2}} \right] \right|^2\, ds \\
			&\leq \Delta^2 h^{2\underline{c}}\log(h)^2 C(\underline{c},\overline{c}),
		\end{split} \label{eqn:technical-1}
	\end{align}
	and 
	\begin{align}
		\int_0^h \left| (h-s)^{a_s-\frac{1}{2}} - (h-s)^{b_s-\frac{1}{2}} \right|^2\, ds 
		&\leq \Delta^2 h^{2\underline{c}}\log(h)^2 C(\underline{c},\overline{c})\label{eqn:technical-2}.
	\end{align}
\end{lemma}
\begin{proof}
	Let $h\in (0,1)$. For some value $c_{sh}$ between $a_{sh}$ and $b_{sh}$, the mean value theorem yields
	\begin{align*}
		&\quad \int_{-\infty}^{0} \left| \left[(h-s)^{a_s-\frac{1}{2}} - (-s)^{a_s-\frac{1}{2}} \right] - \left[(h-s)^{b_s-\frac{1}{2}} - (-s)^{b_s-\frac{1}{2}} \right] \right|^2\, ds \\
		&= \int_{-\infty}^{0}  \left| h^{a_s-\frac{1}{2}} \left[(1-s/h)^{a_s-\frac{1}{2}} - (-s/h)^{a_s-\frac{1}{2}} \right] - h^{b_s-\frac{1}{2}} \left[(1-s/h)^{b_s-\frac{1}{2}} - (-s/h)^{b_s-\frac{1}{2}} \right] \right|^2\, ds \\
		&\hspace{-0.23cm}\overset{\text{subst.}}{=} \int_{-\infty}^{0}  \left| h^{a_{sh}} \left[(1-s)^{a_{sh}-\frac{1}{2}} - (-s)^{a_{sh}-\frac{1}{2}} \right] - h^{b_{sh}} \left[(1-s)^{b_{sh}-\frac{1}{2}} - (-s)^{b_{sh}-\frac{1}{2}} \right] \right|^2\, ds \\
		&= \int_{-\infty}^{0} (a_{sh}-b_{sh})^2\,h^{2c_{sh}} \Big| \log(h)\left[(1-s)^{c_{sh}-\frac{1}{2}} - (-s)^{c_{sh}-\frac{1}{2}} \right] \\
		&\qquad + \left[\log(1-s)\,(1-s)^{c_{sh}-\frac{1}{2}} - \log(-s)\,(-s)^{c_{sh}-\frac{1}{2}} \right] \Big|^2\, ds \\
		&\leq 2\Delta^2h^{2\underline{c}} \log (h)^2 \Bigg[\int_{-\infty}^{0} \left| (1-s)^{c_{sh}-\frac{1}{2}} - (-s)^{c_{sh}-\frac{1}{2}}\right|^2\, ds \\
		&\qquad + \int_{-\infty}^{0} \left| \log(1-s)\,(1-s)^{c_{sh}-\frac{1}{2}} - \log(-s)\,(-s)^{c_{sh}-\frac{1}{2}}\right|^2\, ds \Bigg].
	\end{align*}
	The first integral terms may be bounded as 
	\begin{align*}
		\int_{-\infty}^{0} \left| (1-s)^{c_{sh}-\frac{1}{2}} - (-s)^{c_{sh}-\frac{1}{2}}\right|^2\, ds
		&\leq \int_0^1 \left[ (1+s)^{2\overline{c}-1} + s^{2\underline{c}-1}\right]\, ds + \int_1^\infty s^{2\overline{c}-3}\, ds\\
		&\leq C(\underline{c},\overline{c}),
	\end{align*}
	where we applied the mean value theorem for $s>1$.
	Analogously, 
	\begin{align*}
		&\quad \int_{-\infty}^{0} \left| \log(1-s)\,(1-s)^{c_{sh}-\frac{1}{2}} - \log(-s)\,(-s)^{c_{sh}-\frac{1}{2}}\right|^2\, ds \\
		&\leq \int_0^1 \left[\log(1+s)^2(1+s)^{2\overline{c}-1} + \log(s)^2s^{2\underline{c}-1}\right]\, ds + \int_1^\infty (2+\log(s))^2s^{2\overline{c}-3}\, ds \\
		&\leq C(\underline{c},\overline{c}).
	\end{align*}
	This establishes \eqref{eqn:technical-1}.
	Regarding \eqref{eqn:technical-2}, note that 
	\begin{align*}
		\int_0^h \left| (h-s)^{a_s-\frac{1}{2}} - (h-s)^{b_s-\frac{1}{2}} \right|^2\, ds  
		& \leq \Delta^2 \int_0^h \log(h-s)^2 (h-s)^{2c_s-1}\, ds \\
		& \leq \Delta^2 \int_0^h \log(s)^2 s^{2\underline{c}-1}\, ds \\
		&= \Delta^2 h^{2\underline{c}}\int_0^1 (\log h + \log s)^2 s^{2\underline{c}-1}\, ds.
	\end{align*}
	The latter integral is finite because $\underline{c}\in(0,1)$.
\end{proof}

\begin{proof}[Proof of Theorem \ref{thm:holder-sharp}]
	Theorem \ref{thm:rescaling} implies that almost surely for each $\epsilon>0$ and $\gamma>0$,  \begin{align*}
		\sup_{h\in(0,\epsilon)}\frac{|X_{t+h}-X_t|}{|h|^{H_t+\gamma}} =\infty.
	\end{align*}
	In particular, $\alpha_t(X)\leq H_t$ almost surely, and we conclude that $P(\alpha_t(X)=H_t, \forall t\in[0,T)\cap \Q)=1$.
	Now note that, for all $s$, $t$,
	\begin{align*}
		&\quad \sup \left\{ \alpha\,:\, \lim_{\epsilon\to 0} \sup_{r,v\in (s,t), |r-v|<\epsilon} \frac{|X_r-X_v|}{|r-v|^\alpha} =0\right\} \\
		&\leq \sup \left\{ \alpha\,:\, \lim_{\epsilon\to 0} \sup_{q\in (s,t)\cap \Q} \sup_{0<h<s-q, |h|<\epsilon} \frac{|X_{q+h}-X_q|}{|h|^\alpha} =0\right\} \\
		&\leq \sup \left\{ \alpha\,:\, \sup_{q\in (s,t)\cap \Q} \lim_{\epsilon\to 0} \sup_{0<h<t-q, |h|<\epsilon} \frac{|X_{q+h}-X_q|}{|h|^\alpha} =0\right\} \\
		&\leq \inf_{q\in(s,t)\cap\Q} \alpha_q(X) \quad = \quad \inf_{q\in(s,t)\wedge \Q} H_q,
	\end{align*}	
	which equals $\underline{H}_{[s,t]} = \inf_{r\in[s,t]} H_r$, because $H_t$ is continuous.
	Furthermore, we have shown in \eqref{eqn:Ito-mBm-uniform-Holder} that 
	\begin{align*}
		\sup \left\{ \alpha\,:\, \lim_{\epsilon\to 0} \sup_{r,v\in (s,t), |r-v|<\epsilon} \frac{|X_r-X_v|}{|r-v|^\alpha} =0\right\} \geq \inf_{r\in[s,t]} H_r,
	\end{align*}
	completing the proof.
\end{proof}

\begin{proof}[Proof of Corollary \ref{cor:covariance-local}]
	In the proof of Theorem \ref{thm:rescaling}, all approximations based on the Burkholder-Davis-Gundy inequalities hold not only in probability, but also in $L_2(P)$. 
	That is, for a standard Brownian motion $\tilde{W}_t$ independent of the filtration $\F_t$, we have shown that for some $q>0$ sufficiently small,
	\begin{align*}
		&\quad h^{-H_{t-h^q}} (X_{t+hr}-X_t) \\
		&\deq h^{-H_{t-h^q}} \int_{-\infty}^{\infty} \sigma_{t-h^q} \left[ (t+hr-s)_+^{H_{t-h^q}-\frac{1}{2}} - (t-s)_+^{H_{t-h^q}-\frac{1}{2}} \right]\, d\tilde{W}_s + o_{L_2}(1) \\
		&= Z_{t,r,h} + o_{L_2}(1).
	\end{align*}
	Conditionally on $\F_{t-}$, the process $r\mapsto Z_{t,r,h}$ is a fractional Brownian motion as in \eqref{eqn:def-fbm}, with scaling factor $\sigma_{t-h^q}$ and Hurst parameter $H_{t-h^q}$, the covariance of which is given by \eqref{eqn:cov-fbm}.
	Hence, for $r,v\in\R$,
	\begin{align*}
		\E\left(Z_{t,r,h} \cdot Z_{t,v,h} | \F_{t-}\right) 
		&= \frac{\sigma_{t-h^q}^2 A(H_{t-h^q})}{2} \left( |r|^{2H_{t-h^q}} + |v|^{2H_{t-h^q}} - |r-v|^{2H_{t-h^q}} \right).
	\end{align*}
	Since $\sigma_s$ is continuous and bounded deterministically, we have $\sigma_{t-h^q}\to \sigma_t$ in $L_2(P)$ as $h\to 0$.
	The equicontinuity of $H_t$ yields that $h^{H_{t-h^q}-H_t}\to 1$ and $|r|^{2H_{t-h^q}}\to |r|^{2H_t}$ in $L_2(P)$ as $h\to 0$, as well as $A(H_{t-h^q})\to A(H_t)$ by continuity of $A$.
	Hence, as $h\to 0$,
	\begin{align*}
		&\quad \E\left[ h^{-H_{t-h^q}} \left( X_{t+hr}-X_t \right) h^{-H_{t-h^q}} \left( X_{t+hv}-X_t \right) \right] \\
		&= \E \left[\E\left(Z_{t,r,h} \cdot Z_{t,v,h} | \F_{t-}\right) \right] + o(1)\\
		&\to \E \left[\frac{\sigma_t^2 A(H_t)}{2} \left( |r|^{2H_t} + |v|^{2H_t}-|r-v|^{2H_t} \right) \right].
	\end{align*}
	By the mean value theorem we have
	\begin{align*}
		&\quad  \big|\E\left[ h^{-H_{t}} \left( X_{t+hr}-X_t \right) h^{-H_{t}} \left( X_{t+hv}-X_t \right) \right] \\
		&\qquad - \E\left[ h^{-H_{t-h^q}} \left( X_{t+hr}-X_t \right) h^{-H_{t-h^q}} \left( X_{t+hv}-X_t \right) \right] \big|\\
		&= \left| \E\left[ \left(\exp(-2\log(h){H_{t}})-\exp(-2\log(h){H_{t-h^q}})\right) \left( X_{t+hr}-X_t \right) \left( X_{t+hv}-X_t \right) \right] \right| \\
		&\leq  2 \omega(h^q)\log(h) \exp(2 \omega(h^q)\log(h)) \\
		&\qquad \cdot \sqrt{\E\left[h^{-2H_{t-h^q}}\left| X_{t+hr}-X_t \right|^2 \right]} \sqrt{\E\left[h^{-2H_{t-h^q}}\left| X_{t+hv}-X_t \right|^2 \right]} \\
		&\longrightarrow 0,
	\end{align*}
	as $h \to \infty$. 
	Thus, we obtain the desired result.
\end{proof}

\begin{proof}[Proof of Proposition \ref{prop:cov-stationary}]
	Itô's isometry and the stationarity of $H_t$ and $\sigma_t$ yield
	\begin{align*}
		&\quad \E \left( K_t^H\cdot K_s^H \right)\\
		&= \int_{-\infty}^{\infty} \E\left\{\sigma_r^2 \left[(t-r)_+^{H_r-\frac{1}{2}} - (-r)_+^{H_r-\frac{1}{2}}\right] \left[(s-r)_+^{H_r-\frac{1}{2}} - (-r)_+^{H_r-\frac{1}{2}}\right]\right\} \, dr \\
		&= \int_{-\infty}^{\infty} \E\left\{\sigma_0^2 \left[(t-r)_+^{H_0-\frac{1}{2}} - (-r)_+^{H_0-\frac{1}{2}}\right] \left[(s-r)_+^{H_0-\frac{1}{2}} - (-r)_+^{H_0-\frac{1}{2}}\right]\right\} \, dr \\
		&=  \E \int_{-\infty}^{\infty} \sigma_0^2 \left[(t-r)_+^{H_0-\frac{1}{2}} - (-r)_+^{H_0-\frac{1}{2}}\right] \left[(s-r)_+^{H_0-\frac{1}{2}} - (-r)_+^{H_0-\frac{1}{2}}\right] \, dr \\
		&= \E \left[\frac{\sigma_0^2 A(H_0)}{2} \left( |t|^{2H_0} +|s|^{2H_0}-|t-s|^{2H_0}  \right)\right],
	\end{align*}
	where we applied \eqref{eqn:cov-fbm} in the last step.
	Exchanging the integrals is justified by virtue of Fubini's theorem, since 
	\begin{align*}
		&\quad \int_{-\infty}^{\infty} \sigma_r^2 \left|(t-r)_+^{H_r-\frac{1}{2}} - (-r)_+^{H_r-\frac{1}{2}}\right| \left|(s-r)_+^{H_r-\frac{1}{2}} - (-r)_+^{H_r-\frac{1}{2}}\right| \, dr \\
		&\leq \sqrt{ \E(K_t^H)^2} \sqrt{\E(K_s^H)^2},
	\end{align*}
	which is finite by assumption, using the identity for $\E(K_t^H)^2$ which we have already established.
\end{proof}

\bibliography{mBm}
\bibliographystyle{plain}

\end{document}